\documentclass[a4paper,reqno]{amsart}
\usepackage[english]{babel}
\usepackage[T1]{fontenc} 
\usepackage[utf8]{inputenc}
\allowdisplaybreaks  
\usepackage{amssymb}

\usepackage{faktor}
\usepackage[svgnames,x11names,dvipsnames]{xcolor}
\usepackage{tikz-cd}
\tikzset{isom/.style={draw=none,every to/.append style={edge node={node [sloped, allow upside down, auto=false]{$\simeq$}}}}}

\usepackage{bbm}
\usepackage{mathtools}
\usepackage{needspace}
\usepackage{microtype}
\usepackage[colorlinks,linkcolor=Red3,citecolor=Red3,urlcolor=Cyan3]{hyperref}
\usepackage{orcidlink}

\usepackage{enumitem}
\setlist[itemize,1]{leftmargin=*}
\setlist[enumerate,1]{leftmargin=*,label=\upshape\bfseries\arabic*.}
\setlist[enumerate,2]{leftmargin=*,align=left,label=\upshape\bfseries(\roman*),widest=iii}

\theoremstyle{theorem}
\newtheorem{theorem}[equation]{Theorem}
\newtheorem{corollary}[equation]{Corollary}

\newtheorem*{lemma*}{Lemma}
\newtheorem{proposition}[equation]{Proposition}

\theoremstyle{example}
\newtheorem{definition}[equation]{Definition}
\newtheorem{remark}[equation]{Remark}
\newtheorem{example}[equation]{Example}

\numberwithin{equation}{section} 

\DeclareMathOperator{\Cov}{Cov}
\DeclareMathOperator{\chr}{char}
\DeclareMathOperator{\coker}{coker}
\DeclareMathOperator{\Hom}{Hom}
\DeclareMathOperator{\img}{Im}
\DeclareMathOperator{\Aut}{Aut}
\DeclareMathOperator{\Pic}{Pic}
\DeclareMathOperator{\Har}{\mathbb{H}}
\newcommand{\Het}[1][1]{H^{#1}_{\text{\upshape et}}}
\DeclareMathOperator{\Hrat}{\mathbb{H}_{\text{\upshape ratl}}}
\DeclareMathOperator{\Jac}{Jac}
\DeclareMathOperator{\Div}{Div}
\DeclareMathOperator{\Gal}{Gal}
\DeclareMathOperator{\Ann}{Ann}
\DeclareMathOperator{\Spec}{Spec}
\DeclareMathOperator{\supp}{supp}
\DeclareMathOperator{\Maps}{Maps}
\DeclareMathOperator{\Ram}{Ram}

\renewcommand\emptyset{\varnothing}
\DeclarePairedDelimiterXPP\laurent[1]{}{\lparen\mkern-3.5mu\lparen}{\rparen\mkern-3.5mu\rparen}{}{\ifblank{#1}{\cdot}{#1}}
\DeclarePairedDelimiterXPP\pairing[1]{}{\langle}{\rangle}{}{\ifblank{#1}{\cdot}{#1}}
     
\newcommand{\A}{\mathbb{A}}
\newcommand{\B}{\mathbb{B}}
\newcommand{\C}{\mathbb{C}}
\newcommand{\cc}{\mathcal{C}}
\newcommand{\f}{\mathbbm{f}}
\renewcommand{\H}{\mathbb{H}}
\newcommand{\I}{\mathbb{I}}
\renewcommand{\k}{\Bbbk}
\renewcommand{\ll}{\mathcal{L}}
\newcommand{\lf}{\mathfrak{L}}
\newcommand{\e}{\mathfrak{e}}
\newcommand{\m}{\mathfrak{m}}
\newcommand{\R}{\mathfrak{R}}
\newcommand{\oo}{\mathcal{O}}
\renewcommand{\t}{\mathbbm{t}}
\renewcommand{\P}{\mathbb{P}}
\renewcommand{\ss}{\mathcal{S}}
\newcommand{\Z}{\mathbb{Z}}
\let\stdphi\phi
\let\phi\varphi
\let\varphi\stdphi

\newcommand{\eqdef}{\coloneqq}
\newcommand{\suchthat}{:}
\newcommand{\gkummer}[3][n]{#2\lbrace\mkern0mu{#3}^{1/#1}\mkern0.5mu\rbrace}

\newcommand{\Kxn}[2][n]{K_{x}\lbrace\mkern0mu{#2}^{1/#1}\mkern0.5mu\rbrace}
\newcommand{\AXn}[2][n]{\A_{X}\lbrace\mkern0mu{#2}^{1/#1}\mkern0.5mu\rbrace}

\newcommand{\Zp}{\mathbb{Z}/(p)}
\newcommand{\Zps}{(\mathbb{Z}/(p))^{*}}
\newcommand{\Zns}{(\mathbb{Z}/(n))^{*}}

\newcommand{\joinrelshort}{\mathrel{\mkern-8mu}}
\newcommand{\shortlongrightarrow}{\relbar\joinrelshort\rightarrow}
\newcommand{\isomto}{\mathrel{\mathop{\setbox0\hbox{$\mathsurround0pt\shortlongrightarrow$}\ht0=0.3\ht0\box0}\limits^{\hspace{-1pt}\scalebox{1.0}{$\sim$}\mkern2mu}}}

\newcommand{\correspondence}{\overset{1:1}\leftrightsquigarrow}

\newcommand*{\prodprime}{\operatornamewithlimits{%
\mathchoice
{\prod\nolimits\raisebox{1.618ex}{\hspace{0em}\makebox[0pt]{$\scriptstyle\prime$}}\hspace{-0.2em}}
{\prod\nolimits\raisebox{1ex}{\hspace{-0.15em}\makebox{$\scriptstyle\prime$}}\hspace{0.1em}}
{\prod\nolimits\raisebox{0.618ex}{\hspace{-0.16em}\makebox{$\scriptstyle\prime$}}\hspace{0.1em}}
{\prod\nolimits\raisebox{0.3ex}{\hspace{-0.2em}\makebox{$\scriptstyle\prime$}}\hspace{0em}}
}}

\title[Cyclic covers from an adelic viewpoint]{%
Cyclic covers of an algebraic curve from an adelic viewpoint
}%
\author[L.~M.~Navas~Vicente]{Luis Manuel Navas Vicente\orcidlink{0000-0002-5742-8679}}
\email{navas@usal.es}

\author[F.~J.~PLaza~Mart\'in]{Francisco J. Plaza Mart\'in\orcidlink{0000-0002-5532-7567}}
\email{fplaza@usal.es}

\thanks{Author 1 ORCID: 0000-0002-5742-8679 Author 2 ORCID: 0000-0002-5532-7567.}

\thanks{Research of both authors supported by grant PID2023-150787NB-I00, and of the first author also by grant PID2021-124332NB-C22, both from the Ministerio de Ciencia e Innovación (Government of Spain).}

\thanks{The authors have no relevant financial or non-financial interests to disclose. No datasets are involved in any way in this paper, which only involves research in pure mathematics.}
\address{Departamento de Matem\'aticas and IUFFyM, Universidad de
Salamanca,  Plaza de la Merced 1-4
        \\
        37008 Salamanca. Spain.
}

\begin{document}

\begin{abstract}
We propose an algebraic method for the classification of branched Galois covers of a curve $X$ focused on studying Galois ring extensions of its geometric adele ring $\A_{X}$.  As an application, we deal with cyclic covers; namely, we determine when a given cyclic ring extension of $\A_{X}$ comes from a corresponding cover of curves $Y \to X$, which is reminiscent of a Grunwald-Wang problem, and also determine when two covers yield isomorphic ring extensions, which is known in the literature as an equivalence problem. This completely algebraic method permits us to recover ramification, certain analytic data such as rotation numbers, and enumeration formulas for covers.
\end{abstract}

\maketitle


\section{Introduction}
\label{sec:intro}

The classification of finite extensions of the function field $\Sigma$ of an algebraic curve $X$ is a long standing problem with roots in the foundational work of Riemann, Klein, Hurwitz, and many others. The standard techniques for dealing with this question ultimately rest upon analytic methods in the theory of Riemann surfaces, appearing often in the guise of complex algebraic geometry.

Here we present a purely algebraic alternative approach based on the use of algebraic extensions of the geometric adele ring $\A_{X}$ of an algebraic curve $X$ over an algebraically closed base field $\k$. The motivation is the following fact: a finite Galois cover $Y \to X$ with group $G$ naturally endows $\A_{Y}$ with an $\A_{X}$-algebra structure together with an action of $G$, which makes $\A_{Y}$ a Galois ring extension of $\A_{X}$ with group $G$ in the sense of Chase-Harrison-Rosenberg (CHR)~\cite{ChaseHarrisonRosenberg}.

The main object in our toolkit is the Harrison set 
$\H(R,G)$, defined for any commutative ring $R$ and finite group $G$ (see \S\ref{subsec:galois theory of rings}). It consists of isomorphism classes of Galois ring extensions with fixed Galois group $G$. For abelian groups, $\H(R,G)$ can be itself endowed with a group structure.

Although there is a geometric interpretation of the Harrison set in the language of torsors and the algebraic fundamental group, we choose to work with the Harrison group $\H(\A_{X},G)$ rather than $G$-torsors over $\Spec(\A_{X})$, since describing the full spectrum requires ultrafilters, which unnecessarily complicates matters. In our method we can avoid having to deal with the complete spectrum (\cite[Lemma 3.13, Theorem 3.14]{adeles02}). For the interested reader, the survey~\cite{PaquesGaloisTheories} discusses the connections between the algebraic foundation of CHR and subsequent generalizations or alternative points of view, including the Galois-Grothendieck theory.

Using these algebraic methods, we classify the cyclic covers of $X$ of prime order different from $\chr(\k)$. As must be the case, we recover the bijection between algebraic curves which are finite $G$-Galois covers branched over $n$ points and conjugacy classes of $n$-tuples $(g_{1},\dots,g_{n})$ generating $G$ and satisfying $\prod_{i} g_{i} = 1$. 
The essential tool for this is the Kummer sequence for ring extensions, described for example in~\cite[Theorem 5.4]{Greither} and, from the point of view of étale cohomology in~\cite[Chapter III, Proposition 4.11]{MilneEtaleCohomology}. For readers accustomed to the latter, we can show that $\Har(R,\cc_{n}) \simeq \Het(\Spec R,\cc_{n})$ when $n$ is coprime with $\chr(R)$, thus establishing an equivalence with our approach (but see the remark above regarding $\Spec(\A_{X})$).

The paper is structured around three main themes: the first is determining when, for a prime $p$, given a $\cc_{p}$-Galois ring extension of $\A_{X}$, there exists a $\cc_{p}$-Galois cover of curves $Y \to X$ (\S\ref{subsec:the fundamental exact sequence}) which gives rise to it in the sense explained below. We shall refer to this as the \emph{existence problem}.

The second deals with what is known as an \emph{equivalence problem}, in this case, for function fields of algebraic curves over algebraically closed fields (see~\cite{Sutherland} for a current survey of many instances of this question), namely, when two such covers induce isomorphic extensions of $\A_{X}$ (\S\ref{subsec:the geometric adelic equivalence problem}). 

The third, and final question we will deal with, is how geometric data (such as ramification, analytic rotation numbers, enumeration of covers, etc.) can be recovered from the constructions by completely algebraic means (\S\ref{sec:the role of ramification} and \S\ref{sec:geometric applications}).

The first two problems are solved in \S\ref{sec:Relations between the Kummer Theory of AX and Sigma}, by relating the Kummer theory of the function field $\Sigma$ with that of its adele ring $\A_{X}$. To be precise, the map induced by tensoring with $\A_{X}$ induces a group homomorphism $\Har(\Sigma,\cc_{p}) \to \Har(\A_{X}, \cc_{p})$. The determination of its kernel and image by the exact sequence~\eqref{E:Harrison groups sequence} in Theorem~\ref{T:exact sequence for Harrison} solves the existence problem and also provides the key to answering the equivalence problem, which admits various forms (Corollaries~\ref{C:Harrison groups sequence split},~\ref{C:characterization Hrat equivalences}, and~\ref{C:JacX[p] -> H(Sigma,Cp) -> zero sum}).

Motivated by the fact that the notion of ramification locus appears naturally in the correspondence between finite extensions of the function field $\Sigma$ and covers of the curve $X$, we are led to consider a refinement of~\eqref{E:Harrison groups sequence} in \S\ref{sec:the role of ramification}. For a finite nonempty subset $\R$ of closed points of $X$ representing the ramification locus, we study the Harrison group of the ring corresponding to the affine curve $X \setminus \R$. The main results are encapsulated in the commutative cube~\eqref{E:cube} of Theorem~\ref{T:cube} and applied in \S\ref{subsec:algebraic vs geometric ramification} to obtain various filtrations of Harrison groups by ramification.

The third question is taken up in the last section of the paper (\S\ref{sec:geometric applications}), which presents various situations illustrating how our algebraic approach indeed recovers interesting geometric information. In \S\ref{subsec:classification of covers} we show how our study of the Harrison group corresponds with the classification of $p$-cyclic Galois covers in~\cite{Borevich} and~\cite{Cornalba}. In \S\ref{subsec:enumeration of covers} we recover the enumeration formulas for $p$-cyclic covers with specified ramification given in~\cite{KwakLeeMednykh, Mednykh, Mednykh2}. 
In \S\ref{subsec:rotation data} we show how our methods allow one to define so-called rotation data purely algebraically in terms of the local Kummer symbols at ramified points. We prove that for $\k = \C$, our definition coincides with the usual analytic notion (see for example~\cite{GonzalezDiez}).

We believe that the approach presented here based on studying Galois extensions of the adeles, and in general, of finite $\A_{X}$-algebras, provides an interesting alternative algebraic perspective on  several geometric problems, and, although limited here to cyclic extensions as an illustration,  can be adapted to more general situations.

\subsection{Geometric adeles}
\label{subsec:adeles}

Let $X$ be a projective, irreducible, non-singular curve over an algebraically closed field $\k$. Let $\Sigma$ be the function field of $X$. The notation $x \in X$ will be used to denote a closed point. $A_{x}$, $K_{x}$, $\m_{x}$, $\upsilon_{x}$,  will denote respectively the completion of the valuation ring at $x \in X$, its quotient field, maximal ideal, and the completed valuation. Since $\k$ is algebraically closed, the residue field is $A_{x}/\m_{x} = \k$. 

Recall that the adele ring $\A_{X}$ is the restricted direct product
\begin{equation*}
\label{E:adele}
	       \A_{X} 
	\eqdef \prodprime_{x \in X} (K_{x},A_{x})
	 = \left\{ (\alpha_{x})_{x\in X} \suchthat \alpha_{x}\in A_{x}  \mbox{ for almost all } x\in X \right\}
\end{equation*}
equipped with the restricted product topology. $\A_{X}$ is the direct limit over finite subsets $F$ of (closed) points of $X$ containing an arbitrary fixed $F_{0}$,
\begin{equation}
\label{E:AX limit}
	\A_{X} \simeq \varinjlim_{F \supseteq F_{0}} \A_{X,F},
	\qquad
		   \A_{X,F}
	\eqdef \prod_{x \in F} K_{x}
	\times \prod_{x \in X \setminus F} A_{x}.
\end{equation}
The limit topology is generated by the neighborhood basis at $0$ consisting of the sets $\prod_{x} \m_{x}^{n_{x}}$, where $(n_{x})_{x\in X}$ runs over collections of non-negative integers $n_{x}$ such that $n_{x}=0$ for almost all $x$. It coincides with the restricted product topology. A third characterization is in terms of the $\k$-vector subspaces commensurable with $\A_{X}^+ \eqdef \prod_{x \in X} A_{x}$ as considered by Tate.

The idele group $\I_{X}$ of invertible elements of $\A_{X}$ is the restricted product of $K_{x}^{*}$ with respect to the unit groups $A_{x}^{*}$ and also the direct limit $\varinjlim \I_{X,F}$ where $\I_{X,F} \eqdef \A^{*}_{X,F}$.


\subsection{Galois theory of rings}
\label{subsec:galois theory of rings}

Let us give a brief overview of the Galois theory of commutative ring extensions. Standard sources for this include ~\cite{Borevich,ChaseHarrisonRosenberg,Greither}. A far-reaching survey of this topic and how it has evolved may be found in~\cite{PaquesGaloisTheories}.
We have extracted the following definitions and results from our previous paper~\cite{adeles02} and copied them here for the reader's convenience.

We begin by recalling the definition of a Galois extension of commutative rings (taken from~\cite[Theorem 1.3]{ChaseHarrisonRosenberg}, which generalizes the classical case of fields.

\begin{definition}[Galois extension of rings]
\label{D:Galois extension of rings}
A Galois extension of a commutative ring $R$ consists of a pair $(S,G)$, where $S$ is a commutative ring extension of $R$ and $G$ is a finite group which acts faithfully on $S$ by $R$-algebra automorphisms, with invariants $S^{G} = R$, and satisfying one of the following equivalent conditions (we only list the ones which we will use later on):
\begin{enumerate}

\item $S$ is a separable $R$-algebra and the elements of $G$ are pairwise strongly distinct.

\item The map $h : S \otimes S \to \Maps(G,S)$ given by $h(s \otimes t)(g) = s g(t)$ is an $S$-algebra isomorphism.

\item If $g \in G$ with $g \neq 1$, for any maximal ideal  $\m$ of $S$ there is some $s \in S$ such that $g(s) \not\equiv s \bmod \m$.

\end{enumerate}
In this case we say that $S$ is a Galois ring extension of $R$ with Galois group $G$, or simply a $G$-Galois extension of $R$.
\end{definition}

Recall that the notion of strongly distinct for a pair of morphisms of commutative rings $f,g: S \to T$ means that for every nonzero idempotent $e\in T$  there exist $s\in S$ such that $f(s)e \neq g(s)e$.

We will denote Galois extensions $S$ of the commutative ring $R$ with group $G$ acting via $R$-algebra automorphisms as pairs $(S,G)$.

\begin{definition}[Kummerian ring]
\label{D:Kummerian ring}
Let $n$ be a natural number. A commutative ring $R$ is $n$-\emph{Kummerian} if $n$ is prime to $\chr(R)$ and its unit group $R^{*}$ contains a distinguished $n$-cyclic subgroup $\mu_{n}$.
\end{definition}

A field $K$ with $\chr(K)$ prime to $n$ containing the $n$-th roots of unity $\mu_{n} = \mu_{n}(K^{*})$, is $n$-Kummerian and this is the \emph{only} possible choice of subgroup.

In our context, having fixed an algebraically closed field $\k$, choosing $\mu_{n} \eqdef \mu_{n}(\k^{*})$, the group of $n$th roots of unity in $\k$, induces the structure of an $n$-Kummerian ring on any $\k$-algebra $R$ of characteristic prime to $n$. This is compatible with $\k$-algebra morphisms. 

Note that for the adele ring, $\k \subseteq \Sigma$ is diagonally embedded in $\A_{X}$ and we also have copies of $\k$ in each completion $K_{x}$. This is an example of how there may be infinitely many choices of subgroups $\mu_{n} \subseteq R^{*}$, highlighting the need to specify one.

A particular type of cyclic Galois ring extensions of a Kummerian ring may be constructed as follows.

\begin{definition}[$(G, \chi)$-Kummer extensions]
\label{D:G-chi-Kummer extension}
For a fixed $n$-Kummerian base ring $R$, a \emph{$(G,\chi)$-Kummer extension of $R$} is a triple $(\gkummer{R}{u},G,\chi)$, where 
\[
	\gkummer{R}{u} \eqdef R[T]/(T^{n} - u),
\]
with $u \in R^{*}$, and $G$ is an $n$-cyclic group which acts on $S$ via the character $\chi : G \to \mu_{n} \subseteq R^{*}$ by
\begin{equation}
\label{E:chi-action on T}
		g(T) \eqdef \chi(g) T.
\end{equation}
That this in fact is a Galois extension is shown in~\cite[p.20]{Greither}.
\end{definition}

\begin{definition}[Equivariant Isomorphism]
\label{D:equivariantly isomorphic extensions}
Let $R$ be a commutative ring and $G$ a fixed finite group. For $i=1,2$, let $S_{i}$ be a ring extension of $R$ with a faithful action of $G$ by $R$-automorphisms of $S_{i}$. We say that the pairs $(S_{1},G)$ and  $(S_{2},G)$ are equivariantly isomorphic, or simply $G$-isomorphic, via $\phi$, if $\phi$ is an $R$-algebra isomorphism $\phi: S_{1} \isomto S_{2}$ such that $\phi \circ g = g \circ \phi$ for all $g\in G$. 
\end{definition}

It is clear that a $G$-isomorphism preserves the $G$-Galois property of a ring extension. Harrison showed how to classify the set of  Galois extensions of a given ring $R$ and group $G$ with a fixed action, via the following object.

\begin{definition}[The Harrison group~{\cite[p. 67]{ChaseHarrisonRosenberg}}]
\label{D:Harrison set (R,G)}
Given a group $G$, the set of all $G$-isomorphism classes of $G$-Galois ring extensions $S$ over a fixed base ring $R$ with a fixed faithful action of $G$ is called the \emph{Harrison set} of $(R,G)$ and denoted by $\Har(R,G)$. When $G$ is a finite abelian group, $\Har(R,G)$ can be endowed with a group structure. In this case it is called the \emph{Harrison group}.
\end{definition}

We review some of the basic facts regarding group actions on modules over a ring $R$. Let $R$ be an $n$-Kummerian ring with distinguished subgroup $\mu_{n}$, and $S$ an $R$-module. Suppose $G$ is a finite abelian group of order $n$ acting on $S$ via $R$-module automorphisms. Its dual group $\widehat{G}$ will be identified with $\Hom(G,\mu_{n})$ and its elements referred to simply as characters of $G$. We may consider the decomposition of $S$ with respect to the action of $\widehat{G}$, namely, for $\chi \in \widehat{G}$, we define the $\chi$-eigenspace (or isotypical component)
\begin{equation*}
\label{E:chi eigenspace definition}
	S^{\chi} \eqdef \{ \alpha \in S \suchthat g(\alpha) = \chi(g) \alpha \ \forall g \in G \}.
\end{equation*}
Projection onto the $\chi$-eigenspace is given by $\alpha_{\chi} = e_{\chi} \alpha$, where $e_{\chi}$ is the corresponding idempotent in the group algebra,
\begin{equation*}
\label{E:e-chi idempotent}
	e_{\chi} \eqdef \frac{1}{|G|} \sum_{g \in G} \chi(g^{-1}) g \in R[G].
\end{equation*}
We then have the decomposition
\begin{equation}
\label{E:eigenspace decomposition}
	S = \bigoplus_{\chi\in\widehat{G}} S^{\chi}.
\end{equation}

If $G$ is cyclic of prime order and the action is nontrivial, then each $R$-module $S^{\chi}$ is nontrivial.

With regard to the decomposition~\eqref{E:eigenspace decomposition}, the product of two $G$-Galois extensions $S_{i}/R$ for $i=1,2$ is given by
\begin{equation}
\label{E:Harrison product wrt eigenspaces}
	S_{1} \cdot S_{2} \eqdef \bigoplus_{\chi} (S_{1}^{\chi} \otimes S_{2}^{\chi}),
\end{equation}
where $G$ acts on the summand $S_{1}^{\chi} \otimes S_{2}^{\chi}$ via $g(s_{1} \otimes s_{2}) \eqdef \chi(g)(s_{1} \otimes s_{2})$. One checks that this product factors through equivariant equivalence and thus defines the group law on the Harrison group $\Har(R,G)$.

The neutral element with respect to this product is the so-called trivial $G$-Galois extension, defined by $R^{(G)} \eqdef \bigoplus_{\chi} R$, considered as the set of maps from $G$ to $R$ under the standard sum and product, and with the action of $G$ given by $g((r_{\chi})_{\chi}) = (r_{g^{*}\chi})_{\chi}$ where $g^{*}$ denotes composition with multiplication by $g$.

For simplicity, as in the following result, we will restrict the rank $n$ to be a prime number $p$, different from the characteristic of $\k$.

\begin{definition}[$G$-primitive element]
\label{D:G-primitive}
Let $R$ be a $p$-Kummerian ring and $S$ an $R$-algebra on which a $p$-cyclic group $G$ acts via $R$-automorphisms. Given a \emph{nontrivial} character $\chi : G \to \mu_{p} \subseteq R^{*}$, an element $\alpha \in S$ is called $(G,\chi)$-primitive if:
\begin{itemize}

\item $1, \alpha, \dots, \alpha^{p-1}$ is an $R$-module basis of $S$.

\item $g(\alpha) = \chi(g) \alpha$ for $g \in G$, i.e. $\alpha \in S^{\chi}$.

\end{itemize}
In this case we also say that $\alpha$ is $G$-primitive with character $\chi$.
\end{definition}

The subject of primitive elements in ring extensions is in itself a topic which has been studied by various authors (\cite{Nagahara,Paques}). In this regard, we cite the following result, which is~\cite[Proposition 3.21]{adeles02}.

\begin{proposition}
\label{P:G-primitive equivalences}
Let $S$ be a separable algebra over a $p$-Kummerian ring $R$ on which the $p$-cyclic group $G$ acts via $R$-automorphisms with $S^{G} = R$. Fix a nontrivial character $\chi : G \to \mu_{p} \subseteq R^{*}$. 
For an element $\alpha \in S^{\chi}$, the following are equivalent:
\begin{enumerate}

\item $\alpha$ is $(G,\chi)$-primitive, i.e. $1, \alpha, \alpha^{2}, \dots, \alpha^{p-1}$ is an $R$-module basis of $S$.

\item $\alpha^{p} \in R^{*}$.

\item $\alpha$ is invertible in $S$.

\end{enumerate}
If this is the case, then:
\begin{enumerate}

\setcounter{enumi}{3}

\item $\alpha$ generates $S^{\chi}$ as an $R$-module, and $S$ is a free $R$-module of rank $p$.

\item The characteristic polynomial of $\alpha$ is $C_{\alpha}(T) = T^{p} -\alpha^{p} \in R[T]$. It is separable and generates $\Ann(\alpha)$.

\item $S$ is equivariantly isomorphic to the $(G,\chi)$-Kummer extension (Definition~\ref{D:G-chi-Kummer extension}) $(\gkummer[p]{R}{u},G,\chi)$ for $u = \alpha^{p}$.
\end{enumerate}
\end{proposition}

The basic result in the Kummer theory of ring extensions is the following, which is~\cite[Theorem 5.4]{Greither}.

\begin{proposition}[The general Kummer sequence]
\label{P:Kummer sequence general}
For an $n$-Kummerian ring, $R$, fixing a nontrivial character $\chi : \cc_{n} \to \mu_{n} \subseteq R^{*}$ yields an exact sequence of groups:
\begin{equation}
\label{E:Kummer sequence general}
	1 \xrightarrow{} \faktor{R^{*}}{R^{*n}}
	  \xrightarrow{i_{(R,\chi)}} \Har(R,\cc_{n})
	  \xrightarrow{} \Pic(R)[n]
	  \xrightarrow{} 1.
\end{equation}
Here $i = i_{(R,\chi)}$ sends $u \in R^{*}$ to the $(\cc_{n},\chi)$-Kummer extension $\gkummer[n]{R}{u}$, and the second map is projection onto the $\chi$-eigenspace. This sequence is called the \emph{Kummer sequence}.
\end{proposition}

It follows from the Kummer sequence that if $R$ is $n$-Kummerian and $\Pic(R)$ has trivial $n$-torsion, then
\begin{equation}
\label{E:Kummer sequence isomorphism case}
	\Har(R,\cc_{n}) \simeq \faktor{R^{*}}{R^{*n}}
\end{equation}
and thus every isomorphism class of $\cc_{n}$-Galois extensions of $R$ contains a representative which is $(\cc_{n},\chi)$-Kummer. In particular this holds for the function field $R = \Sigma$ and, more interestingly, also for the geometric adele ring $R = \A_{X}$ (\cite[Theorem 3.11]{adeles02}).

\begin{remark}
\label{R:Aut action Kummer sequence}
As in~\cite[Lemma 3.41]{adeles02}, one sees that $i_{(R,\chi)}$ is equivariant with respect to the action of $\Aut(\cc_{n}) \simeq \Zns$ by 
\begin{equation}
\label{E:Zns action R/Rn}	
	b \in \Zns \longmapsto (u \mapsto u^{b}) \in \Aut(R^{*}/R^{*n})
\end{equation}
and the group isomorphism $\chi^{*} : \Aut(\mu_{n}) \simeq \Zns \isomto \Aut(\cc_{n})$ induced by $\chi : \cc_{n} \isomto \mu_{n} \subseteq R^{*}$.
\end{remark}

\begin{example}
\label{EX:Omegaf}	
Let $p$ be a prime different from $\chr\k$. For a choice of nontrivial character $\chi : \cc_{p} \to \mu_{p}(\k^{*})$, the Kummer sequence yields an isomorphism $\Sigma^{*}/\Sigma^{*p} \isomto \Har(\Sigma,\cc_{p})$ given by sending $f \in \Sigma^{*}$ to the Galois \emph{ring} extension
\begin{equation}
\label{E:Sigma[T]/Tp-f}
	\Omega_{f} \eqdef \Sigma[T]/(T^{p} - f)
\end{equation}
with action via $g(T) = \chi(g) T$, which is a $(\cc_{p},\chi)$-Kummer extension. Thus, without loss of generality, we may work with pairs of the form $(\Omega_{f},\cc_{p})$ when considering general representatives of elements in $\Har(\Sigma,\cc_{p})$.

By Capelli's Theorem, $f \notin \Sigma^{*p}$ if and only if $\Omega_{f}$ is a $p$-Kummer field extension of $\Sigma$. In this case $\cc_{p}$ is identified with the classical Galois group $\Gal(\Omega_{f}/\Sigma)$. Recall that by the Kummer theory of fields, every $p$-cyclic field extension of $\Sigma$ is of this form (see for example Birch's paper in~\cite[Ch. III, \S2]{CaFr}).

If $f \in \Sigma^{*p}$, then $\Omega_{f}$ is the trivial extension $\Omega_{f} \simeq \prod^{p} \Sigma$, and $\cc_{p}$ is identified with a subgroup of the symmetric group $\ss_{p}$ generated by a $p$-cycle, permuting the copies of $\Sigma$.
\end{example}

We shall also need the following construction of cyclic Kummer extensions due to Borevich.

\begin{theorem}[{\cite[\S8, Theorem 2]{Borevich}}]
\label{T:Borevich correspondence n-Kummer}	
For an $n$-Kummerian ring $R$ and fixed nontrivial character $\chi : \cc_{n} \to \mu_{n} \subseteq R^{*}$, there is a 
one-to-one correspondence between $n$-cyclic Kummer extensions $S$ of $R$ and pairs $(\lf,\phi)$, where $\lf \in \Pic(R)[n]$ and $\phi$ is an isomorphism $\phi: \lf^{\otimes n} \isomto R$.
\end{theorem}

The construction sends an $n$-cyclic extension $S/R$ to the pair consisting of the eigenspace $S^{\chi}$ together with the map $\phi : (S^{\chi})^{\otimes n} \to S^{\chi^{n}} = S^{\cc_{p}} = R$, which is easily seen to be an isomorphism. The inverse sends a pair $(\lf,\phi)$ to the quotient of $\bigoplus_{i \geq 0} \lf^{\otimes i}$ by the ideal generated by $\phi(\ell^{\otimes n}) - \ell^{\otimes n}$, where $\ell$ runs over the elements of $\lf$.

In fact this correspondence induces a group isomorphism of $\Har(R,\cc_{p})$ to the set of pairs $(\lf,\phi)$ with group law given by $(\lf_{1},\phi_{1}) \cdot (\lf_{2},\phi_{2}) = (\lf_{1} \otimes \lf_{2}, \phi_{1} \otimes \phi_{2})$ (\cite[\S11, Theorem 1]{Borevich}).

\begin{definition}[Conjugacy of $G$-ring extensions]
\label{D:conjugate extensions}
Let $R$ be a commutative ring. For $i=1,2$, let $S_{i}$ be a ring extension of $R$ with a faithful action of a group $G_{i}$ by $R$-automorphisms of $S_{i}$. We say that $(S_{1},G_{1})$ and  $(S_{2},G_{2})$ are conjugate via $(\phi,\tau)$ if $\phi: S_{1} \isomto S_{2}$ is an $R$-algebra isomorphism and $\tau : G_{1} \isomto G_{2}$ is a group isomorphism, such that $\phi\circ g = \tau(g) \circ \phi$ for all $g\in G_1$. We will denote this relation by $(S_{1}, G_{1}) \sim (S_{2}, G_{2})$.
\end{definition}

Note that the group isomorphism $\tau$ is in fact determined by the $R$-algebra isomorphism $\phi$, namely $\tau(g) = g^{\phi} \eqdef \phi \circ g \circ \phi^{-1}$, although it is convenient to denote it as part of a pair $(\phi,\tau)$ as we are doing here.

\begin{remark}
\label{R:conjugation preserves Galois}
If one looks at the definition, it is immediately clear that being a Galois extension $S/R$ is preserved by conjugation.
\end{remark}

\begin{definition}
\label{D:twist (S,G) aut}
Suppose $(S,G)$ is a ring extension of $R$ with $G$ a group acting faithfully by $R$-automorphisms os $S$. Given an automorphism $\tau \in \Aut(G)$ of $G$, we define the \emph{twist} of $(S,G)$ by $\tau$, denoted by $(S,G)^{\tau}$, to be the same ring extension $S/R$ but with the action of $G$ now given by $g(s) \eqdef \tau(g) s$.
\end{definition}

We will make use of the bifunctoriality of the Harrison group, as expressed in the following result taken from~\cite[Proposition 3.10, Corollary 3.11]{Greither}.

\begin{proposition}
\label{P:Harrison bifunctoriality}
Given a homomorphism of commutative rings $R_{1} \to R_{2}$ and a homomorphism $\tau : G_{1} \to G_{2}$ of finite abelian groups, we have a commutative diagram
\begin{equation}
\label{E:Harrison functorial square}	
	\begin{tikzcd}
	\Har(R_{1}, G_{1}) \ar[r,"\tau^{*}"] \ar[d,"R_{2} \otimes_{R_{1}} -"'] & \Har(R_{1},G_{2}) \ar[d,"R_{2} \otimes_{R_{1}} -"]
	\\
	\Har(R_{2}, G_{1}) \ar[r,"\tau^{*}"'] & \Har(R_{2},G_{2})
	\end{tikzcd}	
\end{equation}
of group homomorphisms of Harrison groups.
\end{proposition}

In particular, conjugacy classes of $n$-cyclic Galois ring extensions correspond bijectively to the quotient set of $\cc_{n}$-isomorphism classes of $\cc_{n}$-Galois extensions modulo $\Aut(\cc_{n})$, i.e.
\begin{equation}
\label{E:correspondence conjugacy p-cyclic mod aut}
	\left\{
	\begin{minipage}[c]{0.35\textwidth}
	\raggedright
	Conjugacy classes of $n$-cyclic Galois extensions $(S,G)$ of $R$
	\end{minipage}
	\right\}
	\correspondence
	\faktor{\Har(R,\cc_{n})}{\Aut(\cc_{n})}.
\end{equation}


\subsection{Covers as Galois ring extensions}
\label{subsec:covers as Galois ring extensions}

Since our main goal is to characterize Galois covers of algebraic curves in terms of Galois ring extensions of the geometric adele ring $\A_{X}$, we need to begin by giving a precise notion of Galois cover of a curve. Note that we include covers which are not necessarily irreducible.

\begin{definition}[Galois cover]
\label{D:Galois cover}
A cover of projective, non-singular algebraic curves, $\pi:Y \to X$, where $X$ is irreducible and $\pi$ is separable, is said to be a Galois cover if the quotient of $Y$ by the action of $\Gal(Y/X)$ is $X$.

Define two covers $Y' \to X$, $Y \to X$ to be equivalent when there is a birational map between them as schemes over $X$ (i.e., there are dense open subschemes $U'$, $U$ of $Y'$, $Y$ respectively and an isomorphism of $U'$ and $U$ over $X$).

Given $X$ as in \S\ref{subsec:adeles}, let $\Cov(X,\cc_{p})$ denote the set of isomorphism classes of $p$-cyclic Galois covers $Y \to X$, where we restrict to $Y$ having connected components which coincide with its irreducible components.
\end{definition}

Note that with this definition, ramification is allowed (in accordance, for example, with the definition of finite Galois branched cover in~\cite[p.125]{Szamuely}). Then each equivalence class has a distinguished representative that satisfies the additional property that its irreducible and connected components coincide.

An \emph{irreducible} cover $Y \to X$, determines a finite separable extension $\Omega/\Sigma$ of function fields. Conversely, given such an $\Omega$, the Zariski-Riemann variety $Y$ of $\Omega$, whose set of closed points is the set of discrete valuations on $\Omega$, determines a cover of $X$. This is the classical equivalence between the category of finitely generated field extensions of $\Bbbk$ of transcendence degree $1$ and the category of nonsingular projective irreducible curves and nonconstant morphisms (e.g. \cite[\href{https://stacks.math.columbia.edu/tag/0BY1}{Theorem 0BY1}]{stacks-project}).

A general cover, which may not be irreducible, determines a finite separable \emph{ring} extension $\Omega$ of $\Sigma$, which is isomorphic to a finite product of field extensions. In this case $Y$ is the disjoint union of the corresponding Zariski-Riemann varieties of these fields and $\Omega$ is the total quotient ring of any affine dense open subscheme of $Y$.

\begin{example}
\label{EX:covers}
Consider $X = \P^{1}$ with homogeneous coordinates $[x_{0},x_{1}]$ and the following two covers. First, $Y = \coprod^{p} X$ with the action of $\cc_{p}$ by cyclic permutation. Second, the cover $Z \subseteq \P^{2}$ given by the zeros of the homogeneous ideal generated by $x_{1}^{p} - x_{2}^{p}$ where $[x_{0},x_{1},x_{2}]$ are homogeneous coordinates on $\P^{2}$, with action $[x_{0},x_{1},x_{2}] \to [x_{0},x_{1},\zeta x_{2}]$ where $\zeta$ is a $p$th root of unity.
Observe that $Z$ is connected and birational to $Y$. The conditions imposed above are satisfied by $Y$ but not by $Z$.
\end{example}

The following is an algebraic version of~\cite[(14.2)]{CaFr}, taking into account the non-irreducible case. The reader may check that the proof (discarding the topological aspects) is analogous to the case of a global field considered there.

\begin{theorem}
\label{T:cassels}
Let $\Omega$ be a finite separable ring extension of $\Sigma$ and let $\A_{Y}$ denote the ring of adeles of $\Omega$. Then
\begin{equation}
\label{E:cassels}
	\A_{X}\otimes_{\Sigma} \Omega  \simeq   \A_{Y}	
\end{equation}
as $\A_{X}$-algebras. Furthermore, $\A_{Y}^{+}$ is the integral closure of $\A_{X}^{+}$ in $\A_{X}\otimes_{\Sigma} \Omega$.
\end{theorem}

A Galois cover $Y \to X$ gives a finite Galois extension of rings $\Omega/\Sigma$ of group $\Gal(Y/X)$. In this case, the isomorphism~\eqref{E:cassels} is equivariant with respect to the action of $\Gal(Y/X)$, and
\begin{equation}
\label{E:AY Gal(Y/X) invariant is AX}
	\A_{Y}^{\Gal(Y/X)} = \A_{X}.
\end{equation}

The following result expresses the basic correspondence between the Kummer theory of the function field $\Sigma$, the Harrison group of isomorphism classes of $\cc_{p}$-Galois ring extensions of $\Sigma$, and covers of the algebraic curve $X$ as we have defined above.

\begin{theorem}
\label{T:Harrison to covers}
Let $X/\k$ be a projective, irreducible, non-singular curve over an algebraically closed field $\k$ with $\chr(\k) \neq p$ and $\Sigma$ its function field. There are canonical identifications
\begin{equation}
\label{E:Harrison to covers}
           \faktor{\big( \Sigma^{*}/\Sigma^{*p}\big) }{\Zps} 
    \simeq \faktor{\Har(\Sigma, \cc_{p})}{\Aut(\cc_{p})} 
    \simeq \Cov(X,\cc_{p}),
\end{equation}
where $\Zps$ acts on $\Sigma^{*}/\Sigma^{*p}$ as in~\eqref{E:Zns action R/Rn}. 
\end{theorem}

\begin{proof}
Recall that throughout we fix a character $\chi : \cc_{p} \to \mu_{p}(\k^{*})$ that serves to determine the action used in defining the Harrison group.

The first correspondence in~\eqref{E:Harrison to covers} is a direct consequence of~\eqref{E:Kummer sequence isomorphism case} and equivariance with respect to the action of $\Aut(\cc_{p})$ (see Remark~\ref{R:Aut action Kummer sequence}). Thus we focus on the second correspondence between extensions and covers.

Given an element of $\Har(\Sigma,\cc_{p})$, i.e. a class of $\cc_{p}$-ring extensions of $\Sigma$, we may choose a representative which is a $(\cc_{p},\chi)$-Kummer extension, namely, of the form $\Omega_{f} = \Sigma[T]/(T^{p} - f)$ as in~\eqref{E:Sigma[T]/Tp-f}, with action via $g(T) = \chi(g) T$. As we discussed above, $\Omega_{f}$ is a finite separable ring, isomorphic to a finite product of field extensions of $\Sigma$.  Let $Y$ be the disjoint union of the Zariski-Riemann varieties of these extensions, endowed with the $\cc_{p}$-action induced by a given isomorphism. This determines a cover $\pi : Y \to X$ as in Definition~\ref{D:Galois cover}. 

Note that $f^{b}$, where $b \in \Zps$, gives an equivalent representative in $\Har(\Sigma,\cc_{p})/\Aut(\cc_{p})$ and does not change $Y$, i.e. it yields the same cover.

Conversely, given a $p$-cyclic Galois cover $\pi:Y \to X$, the localization of the morphism of $\oo_{X}$-algebras $\oo_{X}\to \pi_{*}\oo_{Y}$ at the generic point of $X$ yields a $p$-cyclic Galois ring extension $\Sigma\to \Omega$ and thus a class in $\Har(\Sigma, \cc_{p}) / \Aut(\cc_{p})$.
\end{proof}

\begin{remark}
\label{R:covers trivial element harrison}
A fundamental reason for specifying the conditions on Galois covers as in Definition~\ref{D:Galois cover} is that, under 
the correspondence~\eqref{E:Harrison to covers}, the neutral element of $\H(\Sigma,\cc_{p})$ indeed maps to the trivial cover in that sense (see also Example~\ref{EX:covers}).
\end{remark}

Given a cover $\pi : Y \to X$, the ramification locus $\Ram(\pi)$ is the image of the support of the sheaf of relative differentials $\omega_{Y/X}$. Equivalently, it consists of the points $x \in X$ such that $\# \pi^{-1}(x) < \deg(\pi)$. For a finite subset of closed points $\R \subseteq X$, it is natural to consider
\begin{equation}
\label{E:CovR}
		\Cov_{\R}(X,\cc_{p}) \eqdef \{ \pi \in \Cov(X,\cc_p) \suchthat \Ram(\pi) \subseteq \R \},
\end{equation}
so that 
\begin{equation}
\label{E:stratification Cov CovR}
		\Cov(X,\cc_{p}) = \bigcup_{\R} \Cov_{\R}(X,\cc_{p}).
\end{equation}
We will study the relation between~\eqref{E:stratification Cov CovR} and~\eqref{E:Harrison to covers} in~\S\ref{subsec:algebraic vs geometric ramification}.

\begin{remark}
\label{R:ramified as unramified outside R}
As in the classical topological theory of Riemann surfaces, in the geometric case an element of $\Cov_{\R}(X,\cc_{p})$ corresponds to a member $\pi \in \Cov(X~\setminus~\R, \cc_{p})$ which is unramified, i.e. $\Ram(\pi) = \emptyset$ (see~\cite[Theorem 4.6.4]{Szamuely}).
\end{remark}



\section{Relating the Kummer theories of \texorpdfstring{$\A_{X}$}{AX} and \texorpdfstring{$\Sigma$}{Σ}}
\label{sec:Relations between the Kummer Theory of AX and Sigma}

In~\cite{adeles02} we extensively studied the structure of the Harrison group $\H(\A_{X},\cc_{p})$, classifying the $p$-cyclic Galois extensions of the ring of adeles of the curve $X$. Now, by means of the Harrison group of its function field $\Sigma$, namely $\H(\Sigma,\cc_{p})$, we can now relate this to the classification of $p$-cyclic covers of $X$.

By Proposition~\ref{P:Harrison bifunctoriality}, tensoring $\A_{X} \otimes_{\Sigma} -$ yields a canonical map
\begin{equation}
\label{E:map H(Sigma,Cp) H(AX,Cp)}
    \Har(\Sigma,\cc_{p}) \to \Har(\A_{X},\cc_{p})
\end{equation}
which is equivariant under the action of the automorphism group $\Aut(\cc_{p})$. Our main goal in this section is to study~\eqref{E:map H(Sigma,Cp) H(AX,Cp)}, since it encapsulates the relation between the classical Galois theory of field extensions of $\Sigma$ with that of its adele ring $\A_{X}$.

\subsection{The fundamental exact sequence}
\label{subsec:the fundamental exact sequence}

The \emph{valuation vector} of an idele $u \in \I_{X}$ is defined as
\begin{equation}
\label{E:idele valuation vector}
	\upsilon(u) \eqdef (\upsilon_{x}(u_{x}) \bmod p)_{x} \in  \bigoplus_{x \in X} \Zp.
\end{equation}
Clearly the map $\upsilon$ is a group homomorphism and, since $\k$ is assumed to be algebraically closed, we have
\begin{equation}
\label{E:IXp kernel upsilon}
	\I_{X}^{p} = \ker \upsilon : \I_{X} \to \bigoplus_{x} \Zp
\end{equation}
(see~\cite[(3.29)]{adeles02}).
We summarize several important facts in the following.

\begin{proposition}
\label{P:Harrison Sigma AX square}
There is commutative diagram of groups
\begin{equation}
\label{E:Harrison Sigma AX square}  
\begin{tikzcd}
    \Sigma^{*}/\Sigma^{*p} \ar[r,"f \mapsto \f"] \arrow[d,"i_{(\Sigma,\chi)}"',"\wr"] & \I_{X}/\I_{X}^{p} \arrow[d,"\wr"',"i_{(\A_{X},\chi)}"] \ar[r,"\sim"',"\upsilon"] & \displaystyle\bigoplus_{x \in X} \faktor{\Z}{(p)}
    \\
    \Har(\Sigma,\cc_{p}) \ar[r,"\A_{X} \otimes_{\Sigma} -"] & \Har(\A_{X},\cc_{p})
\end{tikzcd}
\end{equation}
where the bottom arrow is the map~\eqref{E:map H(Sigma,Cp) H(AX,Cp)} which centers our attention, and
the vertical arrows are the first terms in the Kummer sequences~\eqref{E:Kummer sequence general} of $\Sigma$ and $\A_{X}$, respectively. In addition,~\eqref{E:Harrison Sigma AX square} is equivariant with respect to the action of $\Aut(\cc_{p}) \simeq \Zps$ as in~\eqref{E:Zns action R/Rn}
for $R = \Sigma$ or $R = \A_{X}$, and the action on valuation vectors given by multiplication, i.e.
\begin{equation}
\label{E:Zps action valuation vectors}	
	b \in \Zps \longmapsto ((\upsilon_{x})_{x}) \mapsto ((b \upsilon_{x})_{x})
\end{equation}
\end{proposition}

\begin{definition}
\label{D:valuation vector extension}	
Given an extension $(\B,\cc_{p})$ of $\A_{X}$, the \emph{valuation vector} $\upsilon(\B,G,\chi)$ of the triple $(\B,G,\chi)$ is the image of the $\cc_{p}$-isomorphism class of $(\B,\cc_{p})$ under $\upsilon \circ i_{(\A_{X},\chi)}^{-1}$ as in~\eqref{E:Harrison Sigma AX square}. 
\end{definition}

The valuation vector of an extension is an invariant which is explicitly given by
\begin{equation}
\label{E:(G,chi)-valuation vector}
	\upsilon(\B,G,\chi) = \upsilon(\alpha^{p}) = (\upsilon_{x}(\alpha^{p}) \bmod p)_{x} \in \bigoplus_{x \in X} \Zp,
\end{equation}
where $\alpha$ is a $(\cc_{p},\chi)$-primitive element  (Definition~\ref{D:G-primitive} and Proposition~\ref{P:G-primitive equivalences}). The existence of the latter is guaranteed by~\cite[Theorem 3.22]{adeles02}, where it is also shown that the definition does not depend on the choice of $\alpha$.

Note that for a $(\cc_{p},\chi)$-extension $\AXn[p]{\t} = \A_{X}[T]/(T^{p} - \t)$ where $\t \in \I_{X}$, the class of $T$ is a $(\cc_{p},\chi)$-primitive element, and thus
\begin{equation}
\label{E:Cp-chi Kummer valuation vector}
	\upsilon(\AXn[p]{\t},\cc_{p},\chi) = \upsilon(\t).
\end{equation}

In order to proceed with our characterization of $\Hrat(\A_{X},\cc_{p})$, we need an auxiliary result, which is interesting in itself, and points out the role played by the geometry of the curve.

\begin{theorem}[The fundamental exact sequence for $\H(\Sigma,\cc_{p})$ and $\H(\A_{X},\cc_{p})$]
\label{T:exact sequence for Harrison}
Let $X/\k$ be a curve satisfying our initial hypotheses, namely, a projective, irreducible, non-singular curve over an algebraically closed field $\k$ with $\chr(\k) \neq p$. Fix a nontrivial character $\chi : \cc_{p} \to \mu_{p} \subseteq \k^{*}$ and consider the Harrison groups $\Har(\Sigma,\cc_{p})$ and $\Har(\A_{X},\cc_{p})$ of $\cc_{p}$-isomorphism classes of $p$-cyclic Galois ring extensions of the function field $\Sigma$ and the adele ring $\A_{X}$, respectively, where in both cases the $\cc_{p}$-action is via $\chi$. Then there is a canonical exact sequence of groups
\begin{equation}
\label{E:Harrison groups sequence}  
    \begin{tikzcd}[column sep=small]
        0 \rar & \faktor{\I_{X}^{p}\cap\Sigma^{*}}{\Sigma^{*p}}\rar & \Har(\Sigma,\cc_{p})  \rar & 
        \Har(\A_{X},\cc_{p}) \rar & 
        \faktor{\Z}{(p)} \rar & 0.
    \end{tikzcd}
\end{equation} 
\end{theorem}

\begin{proof}
Our assumptions imply that the Picard scheme $\Pic^{0}(X)$ exists and is a projective $\k$-scheme isomorphic to the Jacobian variety $\Jac(X)$.
	
Consider the commutative diagram (of groups)
\begin{equation*}
\begin{tikzcd}
    0 \rar & \Sigma^{*}/\k^{*} \dar{\{\cdot\}^{p}} \rar{d} & \Div(X)
    \dar{\cdot p} \rar & \Pic(X) \dar{\cdot p} \rar & 0
    \\
    0 \rar & \Sigma^{*}/\k^{*}  \rar & \Div(X)
     \rar & \Pic(X)  \rar & 0
 \end{tikzcd}
\end{equation*}
where $\Div(X) = \bigoplus_{x \in X} \Z$ denotes the group of divisors on $X$, the map $d$ sends a function $f \in \Sigma^{*}$ to its divisor, i.e. $\sum_{x} \upsilon_x(f) x$, and the vertical maps are raising to the $p$th power in the first arrow and multiplication by $p$ in the other two. The Snake Lemma yields the exact sequence
\begin{equation}
\label{E:snake div pic p}
    \begin{tikzcd}[column sep=small]
        0 \rar & \Pic(X)[p]   \rar & \faktor{\Sigma^{*}}{\Sigma^{*p}} \rar{\bar{d}} &
        \faktor{\Div(X)}{p(\Div(X))} \rar &
        \faktor{\Pic(X)}{p(\Pic(X))} \rar & 0.
    \end{tikzcd}
\end{equation}
where $\Pic(X)[p]$ denotes the $p$-torsion subgroup of the Picard variety of $X$, which is easily seen to coincide with $\Jac(X)[p]$. The map $\bar{d}$ sends a function $f \in \Sigma^{*}$ to its divisor modulo $p$, i.e, $\sum_{x} \upsilon_x(f) x \bmod p$. By~\eqref{E:IXp kernel upsilon}, we have
\begin{equation}
\label{E:jacobian torsion sequence}
    \begin{tikzcd}[column sep=small]
        \faktor{\Div(X)}{p(\Div(X))}
        \simeq \faktor{\I_{X}}{\I_{X}^{p}}
        \simeq \displaystyle\bigoplus\limits_{x \in X} \faktor{\Z}{(p)},
    \end{tikzcd}
\end{equation}
and therefore
\begin{equation}
\label{E:Jacobian p-torsion isomorphism}
\begin{aligned}
       \Jac(X)[p]
     = \ker(\bar{d})
    &= \big\{
        \bar{f}\in \faktor{\Sigma^{*}}{\Sigma^{*p}} \suchthat \upsilon_x(f)=0 \bmod p \quad \forall x\in X
      \big\}
    \\
    &= \faktor{\I_{X}^{p}\cap\Sigma^{*}}{\Sigma^{*p}}.
\end{aligned}
\end{equation}
Finally, consider the commutative diagram corresponding to the degree map and multiplication by $p$,
\begin{equation*}
\begin{tikzcd}
    0 \rar & \Jac(X) \dar{\cdot p} \rar{i} & \Pic(X)
    \dar{\cdot p} \rar{\deg} & \Z \dar{\cdot p} \rar & 0
    \\
    0 \rar & \Jac(X) \rar{i} & \Pic(X) \rar{\deg} & \Z  \rar & 0
\end{tikzcd}
\end{equation*}
Now, the group $\Jac(X)(\k)$ of points of $\Jac(X)$ with values in $\k$ is divisible~\cite[\href{https://stacks.math.columbia.edu/tag/03RP}{Proposition 03RP}]{stacks-project}, hence the first arrow is surjective and hence we obtain an isomorphism $\Pic(X)/p (\Pic(X)) \simeq \Zp$. Thus using~\eqref{E:Harrison Sigma AX square} we obtain~\eqref{E:Harrison groups sequence}.
\end{proof}

\begin{corollary}
\label{T:Jacobian divisibility}
Given a positive integer $n$ and a divisor $D \in \Div(X)$ with degree equal to a multiple of $n$, there always exists a function $f \in \Sigma^{*}$ such that $D \equiv (f) \bmod n$.
\end{corollary}

\begin{proof}
As we have seen in the proof of the theorem, $\Jac(X)(\k)$ is a divisible group. Thus, given a divisor $D \in \Div(X)$ of degree $d \in n\Z$ and any point $x_{0}$, the composition
\[
\begin{tikzcd}
    \Jac(X)(\k) \arrow[r, "{[n]}"] & \Jac(X)(\k) \arrow[rr, "+D - d x_{0}"] & & \Jac(X)(\k)
\end{tikzcd}
\]
is surjective. Hence there is some $D' \in \Div(X)^{0}$ such that $D - d x_{0} + n D' \sim 0$, where $\sim$ denotes linear equivalence; that is, $D - d x_{0} + n D' = (f)$ for some function $f$ as desired.
\end{proof}
%


\subsection{The geometric adelic equivalence problem}
\label{subsec:the geometric adelic equivalence problem}

In this section we deal with the problem of characterizing when two $p$-cyclic Galois extensions of $\Sigma$ have the same image under the map~\eqref{E:map H(Sigma,Cp) H(AX,Cp)}.

This problem is analogous to the so-called equivalence problem for global fields. The latter has a long history, dating back to the 19th century with Kronecker and Hurwitz. It deals in general with the question of when an invariant associated to a global field might classify it up to isomorphism. For example, one may consider its Dedekind or Weil zeta function, the Galois group, or its adele ring, leading to various related notions of equivalence and counterexamples where non-isomorphic fields turn out to have the same invariant. The first example, concerning the Dedekind zeta function, dates back to Gassmann. Komatsu and Perlis considered the Galois group and the adele ring. The problem constitutes an active area of research to this day, see for example~\cite{Sutherland} for a survey of these classical results and the current state of the art.

\begin{definition}
\label{D:Hrat}
Denote by $\Hrat(\A_{X},\cc_{p})$ the subgroup of $\Har(\A_{X},\cc_{p})$ which is
the image of $\H(\Sigma,\cc_{p})$ under~\eqref{E:map H(Sigma,Cp) H(AX,Cp)}. 
\end{definition}

We may describe $\Hrat(\A_{X},\cc_{p})$ as follows. Given a rational function $f \in \Sigma^{*}$, consider
\begin{equation}
\label{E:AX[T]/Tp-ff}	
	\AXn[p]{\f} = \A_{X}[T]/(T^{p} - \f),
\end{equation}
where $\f = (f_{x})_{x}$ is the idele such that $f_{x} \in K_{x}$ is the germ of $f$ at $x$ and $\cc_{p}$ acts on the class of $T$ via $\chi$. Note that the class of $T$ is a $(\cc_{p},\chi)$-primitive element. By Theorem~\ref{T:cassels}, $\AXn[p]{\f} \simeq \A_{Y}$. This motivates the following definition.

\begin{definition}[rational $(\cc_{p},\chi)$-Kummer extension of $\A_{X}$]
\label{D:Cp-chi-Kummer rational extension}
A $p$-cyclic extension of $\A_{X}$ which is $\cc_{p}$-isomorphic to one of the form~\eqref{E:AX[T]/Tp-ff} will be called a \emph{rational} $(\cc_{p},\chi)$-Kummer extension of $\A_{X}$.
\end{definition}

Thus $\Hrat(\A_{X},\cc_{p})$ consists of the $\cc_{p}$-isomorphism classes of \emph{rational}  $(\cc_{p},\chi)$-extensions of $\A_{X}$. This subgroup can be characterized as a kernel and a cokernel as follows.

\begin{corollary}
\label{C:Harrison groups sequence split}    
The exact sequence~\eqref{E:Harrison groups sequence} splits into the following two short exact sequences:
\begin{equation}
\label{E:Hrat como kernel}
    \begin{tikzcd}[column sep=10pt]
        0 \rar &  \Hrat(\A_{X},\cc_{p}) \rar & 
        \Har(\A_{X},\cc_{p}) \simeq \bigoplus\limits_{x \in X} \faktor{\Z}{(p)} \rar & 
        \faktor{\Z}{(p)} \rar & 0
    \end{tikzcd}
\end{equation} 
\begin{equation}
\label{E:Hrat como cokernel}
    \begin{tikzcd}[column sep=8pt]
        0 \rar & 
        \faktor{\I_{X}^{p}\cap\Sigma^{*}}{\Sigma^{*p}}\rar & 
        \Har(\Sigma,\cc_{p}) \simeq \faktor{\Sigma^{*}}{\Sigma^{*p}} \rar  & \Hrat(\A_{X},\cc_{p}) \rar & 0.
    \end{tikzcd}
\end{equation} 
\end{corollary}

\begin{remark}
\label{R:AutCp action Hrat ker coker}
As we mentioned above, the maps between Harrison groups are equivariant under the action of $\Aut(\cc_{p})$ and hence induce actions on the cokernel in~\eqref{E:Hrat como kernel} and on the kernel in~\eqref{E:Hrat como cokernel}, which is isomorphic to $\Jac(X)[p]$ by~\eqref{E:Jacobian p-torsion isomorphism}.
\end{remark}

The sequence~\eqref{E:Hrat como kernel} determines when a given adelic extension comes from a corresponding extension of the function field. This may be regarded as an adelic analog of the Grunwald-Wang problem, in the form given in~\cite{LorenzRoquette}.

\begin{corollary}
\label{C:characterization Hrat equivalences}
In terms of the valuation vector associated to a $\cc_{p}$-extension $(\B,\cc_{p})$ of $\A_{X}$ (Definition~\ref{D:valuation vector extension}), the previous result translates to the following set of equivalences:
\begin{enumerate}

\item The class of $(\B,\cc_{p})$ belongs to the subgroup $\Hrat(\A_{X},\cc_{p})$.

\item The valuation vector $\upsilon(\B,\cc_{p},\chi) = (v_{x}) \in \bigoplus_{x \in X} \Zp$ satisfies
\begin{equation}
\label{E:zero-sum property}
    \sum_{x} v_{x} \equiv 0 \bmod p.
\end{equation}

\item $\upsilon(\B,\cc_{p},\chi) = \upsilon(f)$ for some function $f \in \Sigma^{*}$. In this case, $(\B,\cc_{p})$ is $\cc_{p}$-isomorphic to $\AXn[p]{\f} \simeq \A_{X} \otimes_{\Sigma} \Omega_{f}$ as in~\eqref{E:Sigma[T]/Tp-f} and~\eqref{E:AX[T]/Tp-ff}.
\end{enumerate}
\end{corollary}

From~\eqref{E:Hrat como cokernel} we see that it is possible to have non-isomorphic field extensions of $\Sigma$ whose adele rings are (topologically)\footnote{We do not discuss topologies on $\A_{X}$-algebras here; see~\cite{adeles02} for details.} isomorphic as $\A_{X}$-algebras. An example is given below. This can be considered the analog for an algebraically closed base field of the known equivalence results for global function fields, which were originally studied by Tate~\cite{Tate-endo} and Turner~\cite{Turner}.

\begin{corollary}[Geometric equivalence problem]
\label{C:JacX[p] -> H(Sigma,Cp) -> zero sum}
There is an exact sequence
\begin{equation}
\label{E:Harrison Sigma zero-sum}
	    0
	\to \Jac(X)[p] 
	\to \Har(\Sigma,\cc_p) 
	\to {\bigoplus_{x \in X}}^{\!\!0} \; \faktor{\Z}{(p)}
	\to 0,
\end{equation}
where the notation indicates zero-sum modulo $p$ tuples.

In particular, when $g > 0$, this shows that the geometric adele ring is not enough in general to classify $\Sigma$-extensions.
\end{corollary}

\begin{proof}
This follows from the exact sequence \eqref{E:Hrat como cokernel}, recalling~\eqref{E:Jacobian p-torsion isomorphism}, and Corollary~\ref{C:characterization Hrat equivalences}.
\end{proof}

We illustrate the negative answer to the geometric equivalence problem with the following explicit example.

\begin{example}
\label{Ex:cyclic covers of elliptic and 2 torsion points}
Let $\k$ be an algebraically closed field of characteristic different from $2, 3$ and set $p=2$. Let $X$ be the normalization of the completion of the affine plane curve of equation
\begin{equation}
\label{E:affine plane curve}
v^2 = (u-1)(u+1)u 
\end{equation}
as a closed subscheme of the affine plane $\Spec \k[u,v]$. A routine computation shows that $X$ is an irreducible, non singular, elliptic curve of equation $x_0 x_2^2 = (x_1-x_0)(x_1+x_0)x_1$ in homogeneous coordinates of the projective plane $\operatorname{Proj}\k[x_0,x_1,x_2]$. Let $\mathbbm{o}\in X$ denote the point $X\cap (x_0)_0$ with coordinates $(0,1,0)$. 
Using the Abel map with base point $\mathbbm{o}$,  the law group of $\Jac(X)$ given by addition of divisors (or tensor product of invertible sheaves) can be transported to $X$ and, in this way, $\mathbbm{o}$ turns out to be the neutral element of the composition law on $X$. The reader may recognize in this construction the addition law in plane cubic defined geometrically by union of points and intersection of lines. Since $X$ has genus $1$, $ \Jac(X)[2]$ consists of $p^{2g} = 4$ points which corresponds via ${\mathcal A}$ to $\mathbbm{o}$ and the three affine points of $X$ whose tangent is parallel to the line $x_1=0$. That is,
\begin{multline*}
\Jac(X)[2] =
\{ \oo_{X}, \oo_{X}\big( (1,-1,0) - (0,1,0)\big) , 
\\
 \oo_{X}\big( (1,0,0) - (0,1,0)\big) ,  \oo_{X}\big( (1,1,0) - (0,1,0)\big) \}.
\end{multline*} 
Now, observe that $\Sigma$, the function field of $X$, is the given by the localization $(\k[u,v]/(v^2-(u-1)(u+1)u))_{(0)}$ so that any function in $\Sigma$ can be obtained as a rational function on $u, v$. Let us consider the rational function $f\eqdef \frac{u-1}{u}\in \Sigma^*$. The divisor of $f$ is given by
\begin{equation*}
(f) = (\frac{u-1}{u}) = (\frac{x_1-x_0}{x_1}) = (x_1-x_0)_0 - (x_1)_0 =
    2 (1,1,0) - 2(1,0,0)
\end{equation*}
so that $\upsilon_x(f) =0 \bmod 2$ for all $x\in X$. Note that $f \notin \Sigma^{*2}$, since if $f = h^{2}$ for some $h \in \Sigma^{*}$ then $h$ would be a rational function with exactly one simple pole and one simple zero, which would imply that the curve $X$ has genus $0$. Proceeding similarly with the function $g\eqdef \frac{u+1}{u-1}$, one concludes that
\begin{equation*}
\faktor{\I_{X}^{2}\cap\Sigma^{*}}{\Sigma^{*2}}
=
\{ 1 , f , g , f g\}. 
\end{equation*}
Hence, Theorem~\ref{T:exact sequence for Harrison} implies that we have the following four non-isomorphic $2$-cyclic Galois extensions of $\Sigma$:
\begin{equation*}
\faktor{\Sigma[T]}{(T^2-1)},
\quad
\faktor{\Sigma[T]}{(T^2-f)},
\quad
\faktor{\Sigma[T]}{(T^2-g)},
\quad
\faktor{\Sigma[T]}{(T^2-fg)}.
\end{equation*}
The latter three are nontrivial and, when tensored by $\A_X$, all yield trivial $2$-cyclic Galois ring extensions of $\A_{X}$. 
Thus we have an example of four distinct elements of $\Har(\Sigma,\cc_{2})$ which map to the identity element in $\Har(\A_{X},\cc_{2})$, in particular, they become isomorphic after tensoring with $\A_{X}$.
\end{example}



\section{The role of ramification}
\label{sec:the role of ramification}

\subsection{The fundamental cube}
\label{subsec:the fundamental cube}

As we outlined in the introduction, we will now consider a finite nonempty subset $\R$ of closed points of $X$, which will represent the ramification locus. By studying the Harrison group of the ring of the affine curve $X \setminus \R$, we obtain a commutative cube~\eqref{E:cube} which refines~\eqref{E:Harrison groups sequence}.

In~\cite[\S2]{adeles02} we introduced a notion of ramification for an adelic algebra $\AXn{\t} = \A_{X}[T]/(T^{n} - \t)$ for an idele $\t \in \I_{X}$ and $n$ coprime to $\chr(\k)$, which we now briefly describe. First, the \emph{ramification locus} of  $\t$ is
\begin{equation}
\label{E:ram(t)}	
    \begin{aligned}
	        \Ram(\t)
	\eqdef& \left\{x\in X \suchthat (n,\upsilon_{x}(t_{x}))\neq n\right\}
	=      \left\{x\in X \suchthat \upsilon_{x}(t_{x}) \not\equiv 0 \bmod n \right\},
    \end{aligned}
\end{equation}
and the corresponding \emph{ramification index} at $x$ is
\[
	e_{x} \eqdef \frac{n}{(n, \upsilon_{x}(t_{x}))}.
\]
The vector $\e=(e_{x})$ of integers is the \emph{ramification profile} of $\t$. Observe that $\Ram(\t) = \{ x \in X \suchthat e_{x} > 1 \}$ is a finite set. By~\cite[Theorem 2.19]{adeles02}, isomorphism classes of $\A_{X}$-algebras of the form $\AXn{\t}$ are classified by their ramification profile, namely, $\AXn{\t_{1}}$ and $\AXn{\t_{2}}$ are isomorphic if and only if $\e_{1} = \e_{2}$.

Recalling the notation in~\S\ref{subsec:adeles}, when $n=p$ is prime it is straightforward to check (see~\cite[Lemma 2.13]{adeles02}) that $x \notin \Ram(\t)$ if and only if the $K_{x}$-algebra $\Kxn[p]{t_{x}}$ is isomorphic to $p$ copies of $K_{x}$. Since the latter is the neutral element of $\H(K_{x},\cc_{p})$, this is also equivalent to $\AXn[p]{\t}$ lying in the kernel of the canonical map $\H(\A_{X},\cc_{p}) \to \H(K_{x}, \cc_{p})$.	

Since every $\cc_{p}$-extension of $\A_{X}$ is isomorphic to a $(\cc_{p},\chi)$-Kummer extension, we have the following stratification by ramification detailed in~\cite[\S3.E, Corollary 3.75]{adeles02}:
\begin{equation}
\label{E:stratification conjugacy p-cyclic ramification}
	\begin{aligned}
	\faktor{\Har(\A_{X},\cc_{p})}{\Aut(\cc_{p})}
	&\correspondence
	\coprod_{\R \subseteq X}
	\left\{
	\begin{minipage}[c]{0.46\textwidth}
	\raggedright
	Conjugacy classes of $p$-cyclic Galois extensions $(\B,G)$ of $\A_{X}$ ramified at $\R$
	\end{minipage}
	\right\}
	\end{aligned}
\end{equation}
where $\R$ ranges over finite subsets of closed points of $X$.

In the following we recall the notation from \S\ref{subsec:adeles}. For a finite subset $\R \subseteq X$ (possibly empty)
\[
		   \A_{X,\R}
	\eqdef \prod_{x \in \R} K_{x} \times \prod_{x \in X \setminus \R} A_{x},
	\quad
	\I_{X,\R} \eqdef \A_{X,\R}^{*}.
\]
In addition, we will also denote
\begin{equation}
\label{E:AR}
	A_{\R} \eqdef H^0(X\setminus \R,\oo_X).	
\end{equation}

\begin{theorem}
\label{T:HarrisonAR}
Let $X/\k$ be a projective, irreducible, non-singular curve over an algebraically closed field $\k$ with $\chr(\k) \neq p$, and $\emptyset \subsetneq \R \subset X$ be a finite nonempty subset (of closed points). The Harrison group $\Har(A_{\R},\cc_{p})$ is characterized as
\begin{equation}
\label{E:HarrisonAR}
	       \Har(A_{\R},\cc_{p})
	\simeq \{f \in \Sigma^{*}/(\Sigma^{*})^{p} : \upsilon_{x}(f) \equiv 0 \bmod p \ \forall x \notin \R\}.
\end{equation}
\end{theorem}

\begin{proof}
Consider the exact sequence of sheaves
\begin{equation}
\label{E:OX*Sigma}
	0 \to \oo_{X}^* \to \Sigma^* \to \Sigma^* / \oo_{X}^*\to 0
\end{equation}
and restrict it to the open subscheme $X\setminus \R$. From its long exact sequence of cohomology one obtains
\[
	0 \to \Sigma^* / A_{\R}^*   \to \operatorname{Div}(X \setminus \R) \to \Pic(X \setminus \R) \to 0.
\]
Observe that, since we assume $\R$ to be nonempty, $X \setminus \R$ is affine, equal to $\Spec(A_{\R})$. This yields the following analog to \eqref{E:snake div pic p}:
\begin{equation}
\label{E:PicAS}
	0\to \Pic(A_{\R})[p] \to \faktor{\Sigma^*}{A_{\R}^*(\Sigma^*)^p} 
	\xrightarrow{\Phi} 
	\faktor{\I_{X}}{\I_{X,\R} \I_{X}^p} = \faktor{\Div(A_{\R})}{p( \Div(A_{\R})) } \to 0 
\end{equation}
On the other hand the Kummer sequences~\eqref{E:Kummer sequence general} yield
\begin{equation}
\label{E:ASSigma}	
	\begin{tikzcd}
	0 \arrow[r]  & \faktor{A_{\R}^*}{(A_{\R}^*)^p}   \arrow[r] \arrow[d, hook] & 
	\Har(A_{\R}, \cc_p) \arrow[r]  \arrow[d,"\Psi"'] & 
	\Pic A_{\R}[p] \arrow[r]  \arrow[d] & 0 
	\\
	0 \arrow[r]  & \faktor{\Sigma^*}{(\Sigma^*)^p}   \arrow[r]  & 
	\Har(\Sigma, \cc_p) \arrow[r]   & 
	0 \arrow[r]   & 0,
	\end{tikzcd}
\end{equation}
where $\Psi$ is induced by $\Sigma \otimes_{A_{\R}}-$ and injectivity of the left vertical arrow is clear. The Snake Lemma yields
\[
	0 \to \ker \Psi \to \Pic A_{\R}[p]  
	\xrightarrow{\delta} \faktor{\Sigma^*}{A_{\R}^* (\Sigma^*)^p} \to 
	\coker \Psi \to 0,
\]
where $\delta$ denotes the connecting morphism. Observe that $\delta$ coincides with the map on the l.h.s. of \eqref{E:PicAS}. From the injectivity of the latter, we conclude that $\ker \Psi = 0$ and 
\begin{equation}
\label{E:coker delta}	
	\coker \delta = \coker \Psi = \img \Phi = \faktor{\I_{X}}{\I_{X,\R} \I_{X}^p}.
\end{equation}
Hence, one has the commutative diagram 
\begin{equation*}
	\begin{tikzcd}
	0 \arrow[r]  & \faktor{A_{\R}^*}{(A_{\R}^*)^p}   \arrow[r] \arrow[d, "\simeq"'] & 
	\Har(A_{\R}, \cc_p) \arrow[r]  \arrow[d, hook, "\Psi"'] & 
	\Pic A_{\R}[p] \arrow[r] \arrow[d, hook, "\delta"'] & 0 
	\\
	0 \arrow[r]  & \faktor{A_{\R}^*}{(A_{\R}^*)^p}   \arrow[r]  & 
	\faktor{\Sigma^*}{(\Sigma^*)^p} \arrow[r]   & 
	\faktor{\Sigma^*}{A_{\R}^* (\Sigma^*)^p}    \arrow[r]   & 0 
	\end{tikzcd}
\end{equation*}
which after another application of the Snake Lemma, using~\eqref{E:coker delta}	yields
\begin{equation}
\label{E:HarrisonAR sequence}
	0\to \Har(A_{\R}, \cc_p)  \xrightarrow{\Psi} \faktor{\Sigma^*}{(\Sigma^*)^p}  \to 
	\faktor{\I_{X}}{\I_{X,\R} \I_X^p} \to 0,
\end{equation}
from which~\eqref{E:HarrisonAR} follows immediately.
\end{proof}

\begin{theorem}
\label{T:cube}
Let $X/\k$ be a curve satisfying our initial hypotheses, namely, a projective, irreducible, non-singular curve over an algebraically closed field $\k$ with $\chr(\k) \neq p$. For each finite nonempty subset $\R \subset X$, we have the following commutative cube:
\begin{equation}
\label{E:cube}
\begin{tikzcd}[row sep=1.2em, column sep=1.2em]
\faktor{A_{\R}^*}{(A_{\R}^*)^p} \arrow[dr, hookrightarrow] \arrow[rr] \arrow[dd, hookrightarrow] & & 
\faktor{\I_{X,\R}}{(\I_{X,\R})^p}   \arrow[dr, "\sim"] \arrow[dd, hookrightarrow] &
\\
& \Har(A_{\R},\cc_p)  \arrow[rr, crossing over]  & & \Har(\A_{X,\R},\cc_p)  \arrow[dd, hookrightarrow]
\\
\faktor{\Sigma^*}{(\Sigma^*)^p}   \arrow[dr,"\sim"] \arrow[rr] & & \faktor{\I_{X}}{\I_{X}^p} \arrow[dr, "\sim"] & 
\\
&\Har(\Sigma,\cc_p)  \arrow[rr] \arrow[from=uu, crossing over, hookrightarrow] & & \Har(\A_{X},\cc_p)  \\
\end{tikzcd}
\end{equation}
where the notation is as above. Moreover, the front face of this cube is a cartesian square:
\begin{equation}
\label{E:HarrisonAR fibered product}
	       \Har(A_{\R},\cc_{p})
	\simeq \Har(\Sigma,\cc_{p}) \underset{\Har(\A_{X},\cc_{p})}{\times} \Har(\A_{X,\R},\cc_{p}),
\end{equation}
Furthermore,~\eqref{E:cube} is equivariant with respect to the action of $\Aut(\cc_{p})$ on the various objects, as described in Proposition~\ref{P:Harrison bifunctoriality} and~\eqref{E:Zns action R/Rn}.
\end{theorem}

\begin{proof}
The diagonal maps come from the Kummer sequence~\eqref{E:Kummer sequence general}, while the horizontal and vertical maps arise from the functoriality of $\Har(-,\cc_p)$ in~\eqref{E:Harrison functorial square} and the following (cartesian) diagram:
\begin{equation}
	\begin{tikzcd}
	A_{\R}  \arrow[r] \arrow[d] & \A_{X,\R} \arrow[d]
	\\
	\Sigma  \arrow[r]  & \A_X. 	
	\end{tikzcd}
\end{equation}
The three diagonal isomorphisms follow from the triviality of the respective Picard groups of $\A_{X,\R}$ and $\A_{X}$, shown in~\cite[Theorem 3.11]{adeles02} (and the field case for $\Sigma$). It is easy to check the injectivity of the vertical maps in back. That of the front right vertical map now follows from the one in back.

Finally, the injectivity of $\Har(A_{\R},\cc_{p}) \to \Har(\Sigma,\cc_{p})$ is part of~\eqref{E:HarrisonAR sequence}, which was the injectivity of the map denoted by $\Psi$ in the proof of Theorem~\ref{T:HarrisonAR}. Combining~\eqref{E:HarrisonAR sequence} with the front and right faces of~\eqref{E:cube} we obtain
\begin{equation*}
	\begin{tikzcd}
	  0 \arrow[r]  
	& \Har(A_{\R}, \cc_p)  \arrow[r, "\Psi"] \arrow[d] 
	& \Har(\Sigma, \cc_p) \arrow[r]  \arrow[d] 
	& \faktor{\I_{X}}{\I_{X,\R} \I_X^p} \arrow[r] \arrow[d, equal] 
	& 0 
	\\
	  0 \arrow[r]  
	& \Har(\A_{X,\R},\cc_{p})  \arrow[r]  
	&\Har(\A_{X},\cc_{p}) \arrow[r]
	& \faktor{\I_{X}}{\I_{X,\R} \I_X^p} \arrow[r]
	& 0 ,
	\end{tikzcd}
\end{equation*}
from which~\eqref{E:HarrisonAR fibered product} follows.
\end{proof}

\begin{corollary}
\label{C:cube limit}
There is an isomorphism from the direct limit over finite nonempty subsets $\R \subset X$ of the top floor of~\eqref{E:cube} to the bottom floor.
\end{corollary}

\begin{proof}
The statement is straightforward for the right face of the cube, and for the left face, it means that
\[
	\Har(\Sigma,\cc_p) \simeq \varinjlim_{\R} \Har(A_{\R},\cc_p),
	\quad
	\faktor{\Sigma^*}{(\Sigma^*)^p} \simeq \varinjlim_{\R} \faktor{A_{\R}^*}{(A_{\R}^*)^p}.
\]
The first follows from~\eqref{E:HarrisonAR sequence}, and the second is then a consequence of this and~\eqref{E:ASSigma}.
\end{proof}


\subsection{Algebraic vs. geometric ramification}
\label{subsec:algebraic vs geometric ramification}

The results of the previous section now let us to take ramification into account in the Harrison groups of $\cc_{p}$-extensions of both $\Sigma$ and $\A_{X}$, allowing us to filter by ramification. We end by establishing the concordance between the algebraic and geometric notions of ramification (Proposition~\ref{P:Harrison to covers with ramification}), which underlies the various geometric applications explored in \S\ref{sec:geometric applications}.

Keeping in mind~\eqref{E:Harrison Sigma AX square}, Definition~\ref{D:valuation vector extension} and~\eqref{E:Cp-chi Kummer valuation vector}, and that any class of $\cc_{p}$-extensions of $\A_{X}$ is represented by a $(\cc_{p},\chi)$-Kummer extension (see~\eqref{E:Kummer sequence isomorphism case}), where the character $\chi$ has been fixed beforehand, the following notions are well-defined.

\begin{definition}[Ramification for $\cc_{p}$-extensions of $\A_{X}$]
\label{D:ramification locus AX extension}
Given a class of $\cc_{p}$-Galois ring extensions $\B \in \Har(\A_{X}, \cc_{p})$, we define the ramification locus of $\B$ and the corresponding valuations and indices by
\begin{equation}
\label{E:ramification locus AX extension}
		\Ram(\B) \eqdef \Ram(\t),
	\quad
	\upsilon_{x}(\B) \eqdef \upsilon_{x}(\t),
	\quad
	e_{x}(\B) \eqdef e_{x},
\end{equation}
where $\t \in \I_{X}$ is such that $\B$ is represented by the $(\cc_{p},\chi)$-Kummer extension $\AXn[p]{\t}$.
\end{definition}

\begin{corollary}[Ramification filtration by subgroups]
\label{C:Harrison ramification filtration AXR}
Ramification provides us with a natural filtration by subgroups:
\begin{equation}
\label{E:Harrison ramification filtration AXR}
\begin{aligned}
	        \Har(\A_{X},\cc_{p})
	&\simeq \varinjlim_{\R} \Har(\A_{X,\R},\cc_{p})
	 \simeq \varinjlim_{\R} \left\{\B \in \Har(\A_{X},\cc_{p}) \suchthat \Ram(\B) \subseteq \R \right\}.
\end{aligned}
\end{equation}
For a given finite subset $\R \subseteq X$,
\begin{equation}
\label{E:Harrison number of extensions Ram in R}
	\#\left\{\B \in \Har(\A_{X},\cc_{p}) \suchthat \Ram(\B) \subseteq \R \right\} = p^{\# \R}.
\end{equation}
\end{corollary}

\begin{proof}
For a finite subset $\R \subset X$, we conclude from the isomorphisms on the right face of~\eqref{E:cube}, that the following diagram is commutative.
\begin{equation}
\label{E:Harrison filtration sequences}	
	\begin{tikzcd}
	  \Har(\A_{X,\R}, \cc_{p}) \arrow[r,"\sim"] \arrow[d, hook]
	& \faktor{\I_{X,\R}}{(\I_{X,\R})^p} \arrow[r,"\sim"] \arrow[d, hook] 
	& \bigoplus_{x \in \R} \faktor{\Z}{(p)} \arrow[d, hook]
	\\
	  \Har(\A_{X}, \cc_{p}) \arrow[r,"\sim"] 
	& \faktor{\I_{X}}{(\I_{X})^p} \arrow[r,"\sim"]
	& \bigoplus_{x \in X} \faktor{\Z}{(p)}
	\end{tikzcd}
\end{equation}
where, by Corollary~\ref{C:cube limit}, the bottom row is the limit of the top row.
Thus, by observing that the following three conditions are equivalent,
\begin{enumerate}

\item $\B \in \Har(\A_{X,\R}, \cc_{p})$,

\item $\Ram(\B) \subseteq \R$,

\item $\upsilon_{x}(\B) \equiv 0 \bmod p$ for $x \notin \R$,

\end{enumerate}
we obtain~\eqref{E:Harrison ramification filtration AXR} and~\eqref{E:Harrison number of extensions Ram in R}.
\end{proof}

\begin{corollary}[Stratification by ramification]
\label{C:Harrison ramification stratification AXR}
There is a natural stratification of sets, indexed by ramification:
\begin{equation}
\label{E:Harrison ramification stratification AXR}
	\begin{aligned}
	        \Har(\A_{X},\cc_{p})
	&\simeq \bigsqcup_{\R} \left\{\B \in \Har(\A_{X},\cc_{p}) \suchthat \Ram(\B) = \R \right\}
	\simeq \bigsqcup_{\R} \bigoplus_{x \in \R} \Bigl(\faktor{\Z}{(p)}\Bigr)^{*}.
	\end{aligned}
\end{equation}
For a given finite subset $\R \subseteq X$,
\begin{equation}
\label{E:Harrison number of extensions Ram = R}
	\#\left\{\B \in \Har(\A_{X},\cc_{p}) \suchthat \Ram(\B) = \R \right\} = (p-1)^{\# \R}.
\end{equation}
\end{corollary}

\begin{proof}
This follows analogously to the above by observing that the following three conditions are also equivalent:
\begin{enumerate}

\item $\B \in \Har(\A_{X,\R}, \cc_{p}) \setminus \bigcup_{\R' \subsetneq \R} \Har(\A_{X,\R'}, \cc_{p})$,

\item $\Ram(\B) = \R$,

\item $\upsilon_{x}(\B) \equiv 0 \bmod p$ if and only if $x \notin \R$.
\qedhere
\end{enumerate}
\end{proof}

In view of~\eqref{E:Harrison filtration sequences}, for a finite (possibly empty) subset $\R \subset X$, we may express the conditions $\Ram(\B) \subseteq \R$ and $\Ram(\B) = \R$ respectively by
\[
	(\upsilon_{x}(\B))_{x \in X} \in \bigoplus_{x \in \R} \Z/(p),
	\qquad
	(\upsilon_{x}(\B))_{x \in X} \in \bigoplus_{x \in \R} (\Z/(p))^{*},
\]	
if no confusion arises.

We will also need to define a corresponding notion of ramification for $\cc_{p}$-Galois ring extensions of the function field $\Sigma$.

\begin{definition}[Ramification for $\cc_{p}$-extensions of $\Sigma$]
\label{D:ramification locus Sigma extension}
Given a class of $\cc_{p}$-Galois ring extensions $\Omega \in \Har(\Sigma, \cc_{p})$, we define the ramification locus $\Ram(\Omega)$ of $\Omega$ via the map~\eqref{E:map H(Sigma,Cp) H(AX,Cp)} as $\Ram(\A_{X} \otimes_{\Sigma} \Omega)$.
\end{definition}

\begin{remark}
\label{R:ramification Omegaf}
For $f \in \Sigma^{*}$, consider, as in~\eqref{E:Sigma[T]/Tp-f}, the $\cc_{p}$-extension of $\Sigma$ given by $\Omega_{f} = \Sigma[T]/(T^{p} - f)$ with $\cc_{p}$-action by $g(T) = \chi(g)T$ where $\chi : \cc_{p} \to \mu_{p}(\k^{*})$ is a fixed character. Then, from the definition and~\eqref{E:Cp-chi Kummer valuation vector}, we have
\begin{equation}
\label{E:ramification Omegaf}
	\begin{aligned}
	   \Ram(\Omega_{f}) 
	&= \Ram(\A_{X} \otimes_{\Sigma} \Omega_{f}) 
	 = \Ram(\AXn[p]{\f})
	\\
	&= \Ram(\f)
	 = \{x \in X \suchthat \upsilon_{x}(f_{x}) \not\equiv 0 \bmod p \},
	\end{aligned}
\end{equation}
where $\f = (f_{x})_{x} \in \I_{X}$ is the idele of germs as in~\eqref{E:AX[T]/Tp-ff}.
\end{remark}

\begin{definition}
\label{D:HR(Sigma,Cp)}
The subset of classes of $\cc_{p}$-Galois extensions of $\Sigma$ with ramification contained in $\R$ will be denoted by
\[
	\Har_{\R}(\Sigma,\cc_{p}) \eqdef \{ \Omega \in \Har(\Sigma,\cc_{p}) \suchthat \Ram(\Omega) \subseteq \R \}.
\]
\end{definition}

With a slight abuse of notation (recalling that the front face of~\eqref{E:cube} is a cartesian square), note that $\Ram(\Omega)$ is the intersection of the finite subsets $\R \subset X$ such that $\Omega \in \Har(A_{\R}, \cc_{p})$.

\begin{proposition}
\label{P:ram Sigma extension H(AR,Cp)}
\begin{equation}
\label{E:ram Sigma extension H(AR,Cp)}
	  \Har_{\R}(\Sigma,\cc_{p})
	= \begin{dcases*}
	  \Har(A_{\R},\cc_{p}), & if $\R$ is nonempty,
	  \\
	  \Jac(X)[p],		    & if $\R = \emptyset$.
	  \end{dcases*}
\end{equation}
\end{proposition}

\begin{proof}
If $\R$ is nonempty, this follows from Corollary~\ref{C:cube limit} and~\eqref{E:HarrisonAR fibered product}, while if
$\R = \emptyset$ then it follows from Corollary~\ref{C:JacX[p] -> H(Sigma,Cp) -> zero sum}.
\end{proof}

\begin{remark}
\label{R:filtration by ramification subgroup HR}
A posteriori, it follows that $\Har_{\R}(\Sigma,\cc_{p})$ is in fact a subgroup of the Harrison group $\Har(\Sigma,\cc_{p})$, and there is a filtration by ramification 
\begin{equation}
\label{E:Harrison ramification filtration Sigma}
\begin{aligned}
	        \Har(\Sigma,\cc_{p})
	&\simeq \varinjlim_{\R} \Har_{\R}(\Sigma,\cc_{p})
	 \simeq \varinjlim_{\R} \left\{\Omega \in \Har(\Sigma,\cc_{p}) \suchthat \Ram(\Omega) \subseteq \R \right\}.
\end{aligned}
\end{equation}
analogous to~\eqref{E:Harrison ramification filtration AXR}.
\end{remark}

\begin{proposition}
\label{P:Harrison Sigma zero-sum R}
For a given finite subset $\R \subseteq X$, we have the following refinement of~\eqref{E:Harrison Sigma zero-sum} filtered by ramification:
\begin{equation}
\label{E:Harrison Sigma zero-sum R}
	    0
	\to \Jac(X)[p] 
	\to \Har_{\R}(\Sigma,\cc_{p})
	\to {\bigoplus_{x \in \R}}^{\!\!0} \; \faktor{\Z}{(p)}
	\to 0.
\end{equation}
\end{proposition}

\begin{proof}
This follows by combining~\eqref{E:Harrison Sigma zero-sum},~\eqref{E:ram Sigma extension H(AR,Cp)} and~\eqref{E:Harrison ramification filtration Sigma}.
\end{proof}

We can now relate the geometric notion of ramification of a cover of an algebraic curve with the corresponding algebraic definition that we have given for extensions of its function field.

\begin{proposition}
\label{P:Harrison to covers with ramification}
If a $p$-cyclic $\Sigma$-ring extension $(\Omega,\cc_{p})$ corresponds to the cover $\pi : Y \to X$ as in~\eqref{E:Harrison to covers}, then
\begin{equation}
\label{E:Ram(Omega)=Ram(pi)}
	\Ram(\Omega) = \Ram(\pi).
\end{equation}
In particular, for a finite subset $\R \subseteq X$,
\begin{equation}
\label{E:HarR(Sigma,Cp) to Cov R}
	\faktor{\Har_{\R}(\Sigma, \cc_p)}{\Aut(\cc_{p})} \isomto 
	\Cov_{\R}(X,\cc_{p})
\end{equation}
with notation as in~\eqref{E:CovR}.
\end{proposition}

\begin{proof}
With notation as in Example~\ref{EX:Omegaf} and the proof of Theorem~\ref{T:Harrison to covers}, it suffices to check that for $f \in \Sigma^{*}$, we have $\Ram(\Omega_{f}) = \Ram(\pi)$. Recall that, by~\eqref{E:ramification Omegaf}, $\Ram(\Omega_{f}) = \{x \suchthat \upsilon_{x}(f_{x}) \not\equiv 0 \bmod p\}$. Noting that in general $(\omega_{Y/X})_{y} \simeq (\omega_{Y/\k})/(\m_{y}^{e-1}\omega_{X/\k})$ where $e$ is the ramification index at $x = \pi(y)$, there are two cases to consider:

\begin{itemize}

\item If $\upsilon_{x}(f_{x}) \not\equiv 0 \bmod p$, then $T^{p} - f_{x}$ is irreducible in $A_{x}[T]$, hence the class of $T$ is a local parameter, provided that $\upsilon_{x}(f_{x}) > 0$. In this case, by uniqueness of extensions of  valuations, we have $\upsilon_{x}(f_{x}) = p \upsilon_{y}(T)$ where $y$ is the unique point in $\pi^{-1}(x)$. Thus $x \in \Ram(\pi)$. If $\upsilon_{x}(f_{x}) < 0$ then the argument is similar with $T^{-1}$ as local parameter.

\item If $\upsilon_{x}(f_{x}) \equiv 0 \bmod p$, and $\upsilon_{x}(f_{x}) \geq 0$, then $T^{p} - f_{x}$ splits completely as $\prod_{i=0}^{p-1} (T - \zeta^{i} f_{x}^{1/p})$ where $\zeta$ is a primitive $p$th root of unity and $f_{x}^{1/p}$ is any choice of $p$th root in $A_{x}$. Now $\pi^{-1}(x) = \{y_{1},\dots,y_{p}\}$ and for each $i$ we have $\upsilon_{x}(f_{x}^{1/p}) = \upsilon_{y_{i}}(T)$. Hence $x \notin \Ram(\pi)$. The case where $\upsilon_{x}(f_{x}) \leq 0$ is similar.

\end{itemize}
Finally,~\eqref{E:HarR(Sigma,Cp) to Cov R} follows from the previous results and~\eqref{E:HarrisonAR fibered product}.
\end{proof}

\begin{corollary}
\label{C:ramified covers existence}
Let $p$ be a prime number. Let $X$ be a smooth, irreducible, non-singular curve over an algebraically closed field $\k$ of characteristic different from $p$. Let $\R \subset X$ be a finite subset of closed points containing at least two elements. Then there exist:
\begin{enumerate}

\item An irreducible $p$-cyclic Galois cover of $X$ whose ramification locus is exactly $\R$. 

\item An irreducible polynomial $T^{p} - f \in \Sigma[T]$ which is ramified exactly at the points of $\R$ and whose Galois group is $\cc_{p}$.

\end{enumerate}  
\end{corollary}

\begin{proof}
This follows from the previous result and Corollary \ref{C:Harrison ramification stratification AXR}. We need $\R$ to have at least two elements so that we can have non-zero valuations which sum to $0$.
\end{proof}



\section{Geometric applications}
\label{sec:geometric applications}

The previous sections, concerned with $\cc_{p}$-Galois ring extensions of the function field $\Sigma$ and the corresponding adele ring $\A_{X}$, allow us to recover some classical geometric results regarding $p$-cyclic Galois covers of an algebraic curve. Since we employ purely algebraic techniques, this is done without any underlying appeal to the analytical and topological theory of Riemann surfaces, as is the case in the standard approach.

\subsection{Classification}
\label{subsec:classification of covers}

We can give a second proof of Proposition~\ref{P:Harrison to covers with ramification}, relating algebraic and geometric ramification, by a direct argument of independent interest, since it shows the concordance of two points of view, namely Borevich's correspondence given in Theorem~\ref{T:Borevich correspondence n-Kummer} (and sketched following its statement) and that of Cornalba in~\cite[Lemma 1]{Cornalba}. 
We refine these constructions to take into account a specified ramification locus $\R$, via~\eqref{E:HarR(Sigma,Cp) to Cov R}, with notation as in \S\ref{subsec:algebraic vs geometric ramification}. 

\begin{proposition}[Proposition~\ref{P:Harrison to covers with ramification} via Borevich and Cornalba]
\label{P:Harrison to covers with ramification Borevich Cornalba}
For a finite subset $\R \subseteq X$,
\begin{equation}
\label{E:HarR(Sigma,Cp) to Cov R Borevich Cornalba}
	\faktor{\Har_{\R}(\Sigma, \cc_p)}{\Aut(\cc_{p})} \isomto 
	\Cov_{\R}(X,\cc_{p}).
\end{equation}
\end{proposition}

\begin{proof}
Cornalba~\cite{Cornalba} showed that there is a bijection of sets
\begin{equation}
\label{E:triple}
    \Cov(X,\cc_{p}) 
    \correspondence
    \left\{
    \begin{minipage}[c]{0.41\textwidth}
    \raggedright
    $(\ll,D)$ where $D \in \Div(X)$ is an effective divisor and $\ll \in \Pic(X)$ is a line bundle with $\ll^p \simeq \oo_X(D)$
    \end{minipage}
    \right\}
    \!\!\Biggm/\sim,
\end{equation}
where two pairs are defined to be equivalent, denoted by $(\ll_{1},D_{1})\sim (\ll_{2},D_{2})$ with $D_i= \sum_{x\in X} a^{i}_{x} x$, if there exists  $1\leq b< p$ such that
\begin{itemize}

\item $b a^1_{x} \equiv a^2_{x} \bmod p$ for all $x\in X$, 

\item $\ll_1^b \simeq \ll_2 (\sum c_{x} x)$ with $b a^1_{x} = c_{x} p +  a^2_{x}$ for each $x\in X$.

\end{itemize}
In this correspondence, $Y$ is not required to be irreducible. Indeed, the trivial cover $Y = \coprod^{p} X$ corresponds to the equivalence class of the pair $(\oo_X,0)$. In addition, a given cover $\pi \in \Cov(X,\cc_{p})$ corresponding to a pair $(\ll,D)$ has ramification locus $\Ram(\pi)$ equal to $\supp(D)$.

Fix a character $\chi : \cc_{p} \to \mu_{p}(\k^{*})$, that determines the actions used in defining the different Harrison groups involved, as well as a finite subset $\R \subseteq X$.

Let $\Omega$ be a $\cc_{p}$-ring extension of $\Sigma$ with $\R \supset \Ram(\Omega)$, i.e. whose class lies in $\Har_{\R}(\Sigma,\cc_{p})$. By~\eqref{E:ram Sigma extension H(AR,Cp)} and~\eqref{E:Harrison ramification filtration Sigma}, $\R$ may be assumed nonempty and $\Omega \in \Har(A_{\R},\cc_{p})$. The unramified case, where $\R = \emptyset$, is easily seen to follow from the general result.

Under Borevich's correspondence, there is a pair $(\ll,\phi)$ associated to $\Omega$, where $\ll \in \Pic(A_{\R})[p]$ and $\phi: \ll^{\otimes p} \isomto A_{\R}$. The exact sequence
\[
	\bigoplus_{x_i\in \R} \Z x_{i} \to \Pic(X) \to \Pic(X\setminus \R) \to 0
\]
shows that  there exists $\lf \in \Pic(X)$ such that $\lf \vert_{X\setminus \R} \simeq \ll$. Then, applying Lemma 30.10.6 \cite[\href{https://stacks.math.columbia.edu/tag/0FD0}{Lemma 0FD0}]{stacks-project} to $\phi^{-1}:\oo_{X\setminus \R}\isomto \ll^{\otimes (-p)}$ shows that there exists an effective divisor $\bar D$ with support contained in $\R$ and a morphism of $\oo_X$-modules 
\[
	\bar\phi^{-1} : \oo_X(-\bar D) \to \lf^{\otimes (-p)}
\]
whose restriction to $X\setminus \R$ coincides with $\phi^{-1}$. Choosing a suitable $\bar D$, it can be assumed that $\bar\phi^{-1}$ is an isomorphism. Write $\bar D=p E + D$ for divisors $E,D$ with support contained in $\R$ and such that $D=\sum_{x_i\in \R} a_i x_i$ with $0\leq a_i\leq p-1$. Then, $\bar\phi$ yields
\[
	( \lf(-E) )^{\otimes p} \simeq \oo_X(D).
\]
It is straightforward to check that the pair $(\lf(-E), D )$ is the data associated by Cornalba to the cover corresponding to the extension $\Omega$ and that its branch locus is $\{x_i\in \R \suchthat a_i\neq 0\} = \supp(D) \subseteq \R$. Thus we obtain an element of $\Cov_{\R}(X,\cc_{p})$. If one chooses another character, the resulting pair is equivalent in the sense of \eqref{E:triple}.

Conversely, let $(\ll,D)$ be associated to an element of $\Cov_{\R}(X,\cc_{p})$ as in \eqref{E:triple}; thus $\R \supseteq \supp(D)$. Fix an isomorphism $\phi: \ll^{\otimes p}\isomto \oo_X(D)$, which is unique up to $\Bbbk^*$. Restricting to the open subset $X\setminus \R$, let $\lf \eqdef H^{0}(X\setminus \R, \ll) \in \Pic(A_{\R})$ and $\phi_{\R} \eqdef \phi|_{X \setminus \R} : \lf^{\otimes p} \isomto A_{\R}$. Since $X\setminus \R$ is affine with ring $A_{\R}$, we obtain a Borevich pair $(\lf,\phi_{\R})$, which then corresponds (via the choice of character $\chi$) to an element of $\Har(A_{\R},\cc_{p}) = \Har_{\R}(\Sigma,\cc_{p})$.

Note that a different choice of character yields a conjugate extension, not necessarily an isomorphic one. Thus we obtain an element of $\Har_{\R}(\Sigma,\cc_{p})/\Aut(\cc_{p})$.
\end{proof}

Theorem~\ref{T:Harrison to covers} gave a geometric description of the Harrison group $\H(\Sigma,\cc_{p})$ and now, with the previous discussion in mind, we obtain the following geometrical interpretation of its image $\Hrat(\A_{X},\cc_{p})$ inside $\H(\A_{X},\cc_{p})$.

\begin{theorem}
\label{T:Cornalba Jac[p]}
There is a canonical bijection of sets
\begin{equation}
    \label{E:Cornalba Jac[p]}
    \faktor{\Cov(X,\cc_{p})}{\Jac(X)[p]} 
    \correspondence
    \faktor{\Hrat(\A_{X},\cc_{p}) }{\Aut(\cc_{p})}.
\end{equation}
\end{theorem}

\begin{proof}
First observe that $\Jac(X)[p]$ acts freely on the r.h.s. of~\eqref{E:triple} since $L\in \Jac(X)[p]$ acts by sending a pair $(\ll,D)$ to $(D,\ll\otimes L)$ and, thus, we conclude that 
\begin{equation*}
    \faktor{\Cov(X,\cc_{p})}{\Jac(X)[p]} 
    \correspondence
        \left\{
    \begin{minipage}[c]{0.27\textwidth}
        \raggedright
        $D\in \Div(X)$ such that  $\deg(D)  \equiv 0 \bmod p$
    \end{minipage}
    \right\}
    \!\!\Bigm/\sim,
\end{equation*}
where $D_1= \sum_{x} a^{1}_{x} x $ and  $D_2=\sum_{x} a^2_{x} x$ are equivalent if there exists $b\in \Zps$ such that $b a^1_{x}\equiv a^2_{x}\bmod p$ for all $x$. 

On the other hand, recalling~\eqref{E:jacobian torsion sequence} and~\eqref{E:Hrat como kernel}, it is clear that 
\begin{equation*}
    \label{E:Cornalba}
    \Hrat(\A_{X},\cc_{p}) 
    \simeq
        \left\{
    \begin{minipage}[c]{0.42\textwidth}
        \raggedright
        $D\in \Div(X)$ with coefficients in $\Zp$ such that  $\deg(D)  \equiv 0 \bmod p$
    \end{minipage}
    \right\}.
\end{equation*}
Finally, from~\eqref{E:Jacobian p-torsion isomorphism} one has an isomorphism $ \Jac(X)[p]\simeq \I_{X}^{p}\cap\Sigma^{*} / \Sigma^{*p}$, and now the statement is deduced from the exact sequence~\eqref{E:Hrat como cokernel}, recalling Remark~\ref{R:AutCp action Hrat ker coker}.
\end{proof}


\subsection{Enumeration}
\label{subsec:enumeration of covers}

\begin{theorem}
\label{T:enumeration covers ram R}
Let $X/\k$ be a projective, irreducible, non-singular curve of genus $g$ over an algebraically closed field $\k$ with $\chr(\k) \neq p$. For a finite subset $\R \subseteq X$ with $r$ points,
\begin{enumerate}

\item The number of nontrivial unramified $p$-cyclic covers $\pi \in \Cov(X,\cc_{p})$ is
\begin{equation}
\label{E:number nontrivial unramified covers}
	\frac{p^{2g} - 1}{p - 1}.
\end{equation}

\item The number of nontrivial $p$-cyclic covers $\pi \in \Cov(X,\cc_{p})$ with ramification contained in $\R$, assuming $\R$ nonempty, is

\begin{equation}
\label{E:number nontrivial covers ram <= R}
	\frac{p^{2g+r-1} - 1}{p - 1}.
\end{equation}

\item The number of nontrivial $p$-cyclic covers $\pi \in \Cov(X,\cc_{p})$ with ramification equal to $\R$, assuming $\R$ nonempty, is
\begin{equation}
\label{E:number nontrivial covers ram=R}
	p^{2g-1}( (p-1)^{r-1} + (-1)^{r}).
\end{equation}

\end{enumerate}
\end{theorem}

\begin{proof}
\leavevmode	
\begin{enumerate}

\item By~\eqref{E:ram Sigma extension H(AR,Cp)}, we have $\H_{\emptyset}(\Sigma,\cc_{p}) = \Jac(X)[p]$, which has $p^{2g}$ elements. Note that the isotropy subgroup of a class in $\H_{\emptyset}(\Sigma,\cc_{p})$ under the action of $\Aut(\cc_{p}) \isomto \Zps$ is either trivial or all of $\Zps$. Observing that the former only holds for the trivial cover and recalling~\eqref{E:Harrison to covers} yields~\eqref{E:number nontrivial unramified covers}.

\item Combining~\eqref{E:Harrison Sigma zero-sum}, Corollary~\ref{C:Harrison ramification filtration AXR}, and~\eqref{E:Ram(Omega)=Ram(pi)}, $\H_{\R}(\Sigma,\cc_{p})$ has $p^{2g+r-1}$ elements. Reasoning as before, we obtain~\eqref{E:number nontrivial covers ram <= R}.

\item Considering again~\eqref{E:Harrison Sigma zero-sum} and the stratification by ramification~\eqref{E:Harrison ramification stratification AXR}, the number of elements in $\H(\Sigma,\cc_{p})$ with ramification exactly equal to $\R$ is $p^{2g}$ times the number of solutions to the congruence
\[
	\upsilon_{1} + \dots + \upsilon_{r} \equiv 0 \bmod p,
	\qquad
	\upsilon_{1},\dots,\upsilon_{r} \in \Zps,
\]
which is given by
\begin{equation}
\label{E:number solutions x1+...+xr=0(p)}
	\frac{1}{p}( (p-1)^{r} + (-1)^{r}(p-1) ).
\end{equation}
To see this, for a subset $A \subseteq \R$, let $h(A)$ be the number of tuples $\upsilon : \R \to \Zp$  with support equal to $A$, i.e. $\upsilon(x) \not\equiv 0 \bmod p$ if and only if $x \in A$, and satisfying $\sum_{x \in \R} \upsilon(x) \equiv 0 \bmod p$. We want to compute $h(\R)$. Set $H(A) = \sum_{B \subseteq A} h(A)$, which is the number of zero-sum tuples with support contained in $A$. This is easily seen to equal $p^{\#A - 1}$ when $A$ is nonempty and $1$ when $A = \emptyset$. Note that this only depends on the cardinality of $A$. By Möbius inversion, $h(A) = \sum_{B \subseteq A} (-1)^{\#A - \#B} H(B)$. Rearranging by cardinality we obtain $h(\R) = (-1)^{r} + \sum_{k=1}^{r} \binom{r}{k} (-1)^{r-k} p^{k-1}$, which simplifies to the given formula.

In this enumeration we are not counting the class of the trivial extension, since it is unramified and we assume $\R$ is nonempty. Noting that $\Aut(\cc_{p})$ acts freely on every class with prescribed nonempty ramification $\R$ yields~\eqref{E:number nontrivial covers ram=R} after dividing by $p-1$.
\qedhere
\end{enumerate}
\end{proof}

Alternatively, \eqref{E:number nontrivial covers ram=R} may be deduced from~\eqref{E:number nontrivial unramified covers} and~\eqref{E:number nontrivial covers ram <= R} by an analogous process of Möbius inversion, which leads to
\[
	h(\R) = (-1)^{r} \frac{p^{2g} - 1}{p - 1}
	        + \sum_{k=1}^{r} \binom{r}{k} (-1)^{r-k} \frac{p^{2g+k-1} - 1}{p - 1},
\]
and simplifying. In either case~\eqref{E:Harrison Sigma zero-sum} is the starting point.

Note that~\eqref{E:number nontrivial unramified covers} is~\cite[Theorem 3]{Mednykh} and is also found on~\cite[p.490]{Jones} and the first formula of~\cite[Theorem 7]{KwakLeeMednykh}, while~\eqref{E:number nontrivial covers ram=R} is on~\cite[p. 500]{Jones} and is the second formula in~\cite[Theorem 7]{KwakLeeMednykh}. These papers use different methods, although they ultimately rely on the properties of the fundamental group of a Riemann surface, as is common in this kind of problem. Thus we think it is interesting to point out that our approach, avoiding the dependence on the classical topological or analytic structure, allows us to generalize several aspects of $p$-cyclic Galois covers, such as classification and enumeration, to any algebraically closed base field of any characteristic (different from $p$) via the Galois theory of the adele ring. We may also deal with the general abelian case using these tools.

The reader may compare our methods to the references above, for example~\eqref{E:number solutions x1+...+xr=0(p)} is~\cite[Lemma 3]{Mednykh2} although we have arrived at this congruence via quite a different path.

\begin{corollary}
\label{C:existence of p-covers}
There exists a $p$-cyclic Galois cover of $X$ with exactly $r \geq 0$ ramification points, except in the following two cases: $g = 0$ and $r = 0$, or $g \geq 0$ and $r = 1$.
\end{corollary}

\begin{proof}
This is immediate from the formulas in Theorem~\ref{T:enumeration covers ram R}.
\end{proof}

Recall that the Hurwitz existence problem is the question of whether a certain set of numerical data associated to a curve $X$, foremost among these the Riemann-Hurwitz relation, is in fact realizable by a cover $Y \to X$. This problem has been extensively studied, beginning with Hurwitz himself. It is known that for genus $g > 0$ the answer is affirmative (major advances were made in~\cite{EdmondsKulkarniStong} and many remaining cases dealt with later on), while the case $g = 0$ reduces to $3$-point branched covers of the Riemann sphere and remains a difficult question. In particular when $g > 0$ and $r = 1$, although as we see above, there are no $p$-cyclic Galois covers, other covers do exist as predicted by these general results. For $g = 0$ and $r = 2$,~\eqref{E:number nontrivial covers ram=R} says there is a unique $p$-cyclic cover. For $r=3$ there are $p-2$ such covers.


\subsection{Rotation data}
\label{subsec:rotation data}

In this section we give an algebraic definition of rotation numbers (Definition~\ref{D:algebraic rotation data}), and using the Kummer pairing at ramification points, we prove that for $\k = \C$ it coincides with the classical one for the case of a compact connected Riemann surface. 

To achieve this, we begin by considering the following classical construction of a $p$-cyclic cover of the projective line, which is standard in the literature (see e.g.~\cite[\S5]{BSW} and especially~\cite[\S1]{GonzalezDiez}).

Set $\k= \C$, $X= \P^{1}$, $p$ a prime number, and $\R \eqdef \{x_{1}, \ldots, x_{r}\}$ a finite set of distinct points of $X$. Let $Y$ be the normalization of 
\begin{equation}
\label{E:superelliptic}
    y^{p}  =  f(x) \eqdef ( x- x_1)^{v_1} \cdots (x-x_r)^{v_r},	
\end{equation}
where $\pi:Y\to X$ maps $(x,y)$ to $x$, and the exponents $v_{i}$ satisfy $0 < v_{i} < p$ and $\sum_{i} v_{i} \equiv 0 \bmod p$ (the latter is equivalent to being unramified at $\infty$). Observe that the Riemann-Hurwitz formula implies that $r \geq 1$ in the above expression for $f(x)$.

Fix a nontrivial character $\chi : \cc_{p} \to  \mu_{p} \subseteq \C^{*}$ and let $\zeta = \chi(1) \in \mu_{p}$, which is a primitive $p$th root of unity. This defines an action of $\cc_{p}$ on $Y$ where $1$ acts via the automorphism
\begin{equation}
\label{E:automorphism order p superelliptic}
	\tau : (x,y) \mapsto (x,\zeta y).	
\end{equation}
Thus $\tau$ is a generator of $\Gal(Y/X)$, and allows us to identify $\cc_{p}$ with $\Gal(Y/X)$.

f we assume that  $\operatorname{gcd}(v_1,\ldots, v_r)=1$, then $Y$ is irreducible and $\oo_X\to \pi_*\oo_Y$ at the generic point is
\[
    \Sigma \eqdef \C(x) \hookrightarrow  \Omega_{f} = \Sigma[y]/(y^p-f(x))
\]
as considered in Example~\ref{EX:Omegaf}. This is what we termed a $(\cc_{p},\chi)$-Kummer extension (Definition~\ref{D:G-chi-Kummer extension}). 

Thus $\pi: Y\to X$ is a $p$-cyclic cover with ramification locus
\[
	\Ram(\pi) = \R = \{x_1, \ldots, x_r\},
\]
having ramification index at $x_i$ equal to ${p}/{(p,v_i)}=p$. Thus $\pi$ defines a class in $\Cov_{\R}(X,\cc_{p})$ and
the class of $\Omega_{f}$ is an element of $\H_{\R}(\Sigma,\cc_{p})$, related via the correspondence in Proposition~\ref{P:Harrison to covers with ramification}.

Let us now change our point of view to the geometric adele ring. Under the canonical map~\eqref{E:map H(Sigma,Cp) H(AX,Cp)}, the corresponding adelic algebra $\A_{Y} = \AXn[p]{\f}$ is a rational $(\cc_{p},\chi)$-Kummer adelic extension of $\A_{X}$ (Definition~\ref{D:Cp-chi-Kummer rational extension}), hence lying by definition in $\Hrat(\A_{X},\cc_{p})$. The characteristic polynomial of $y$ is equal to $C_y(T)=T^p - \f$, and thus the associated valuation vector of this extension (Definition~\ref{D:valuation vector extension}) is
\[
      (v_{x})_{x}
	\eqdef \upsilon(\A_{Y},\cc_{p},\chi)
    = (\upsilon_{x}(\f))_{x}
    = \begin{cases}
      v_{i} \bmod p, & x = x_{i} \in \R,
            \\
      0     \bmod p, & x \notin \R,
      \end{cases}
\]
where we take as $\chi$-primitive element the class of $T$ modulo $T^{p} - \f$. This invariant satisfies
\[
    \sum_{x} v_{x} \equiv 0 \bmod p	
\]
in accordance with Corollary~\ref{C:characterization Hrat equivalences}. Recall that these invariants classify $(\cc_{p},\chi)$ adelic extensions and hence $p$-cyclic covers.

In the classical theory of Riemann surfaces, one considers the so-called \emph{rotation numbers}. These are defined in terms of the automorphism $\tau$ given by the action~\eqref{E:automorphism order p superelliptic}, namely, choosing $\zeta = e^{2\pi i/p}$, in a neighborhood of a branch point $x = (x_{i},0)$ one sees that $\tau$ rotates a small disk around this point by an angle $2\pi \rho_{x}/p$; the integer $\rho_{x} \bmod p$ is the corresponding rotation number. A simple calculation shows they are the inverses modulo $p$ of the valuations $\upsilon_{x}$: 
\begin{equation}
\label{E:rotation is inverse of valuation}
	\rho_{x} v_{x} \equiv 1 \bmod p.	
\end{equation}

The usual approach based on complex analytic and topological methods (for example~\cite{GonzalezDiez}), classifies the $p$-cyclic branched covers of the punctured Riemann sphere in terms of the rotation data $(\rho_{x})_{x \in \R}$.

On the other hand, in \S\ref{subsec:the geometric adelic equivalence problem}, we have already classified $p$-covers in arbitrary characteristic prime to $p$ in a purely algebraic manner, reflected explicitly in the valuation vector $(\upsilon_{x})_{x}$. To complete the algebraic picture, we will show how both the valuations and the rotation data arise naturally from the local Kummer symbols at each ramified point $x \in \R$. For this we need to recall some notions introduced in~\cite{adeles02}.

For $x \in \R$, by Kummer theory, there is a unique (up to isomorphism) abelian extension $E_{x}$ of $K_{x}$ of exponent $p$, generated by adjoining the $p$th root of any non $p$th power in $K_{x}$. This is the algebraic analog of the classical Puiseux series expansion in the theory of Riemann surfaces, since for example $K_{x} \simeq \k\laurent{z_{x}}$ where $z_{x}$ is a uniformizer at $x$, and the local valuation $\upsilon_{x}$
\begin{equation}
\label{E:Kummer Kx/Kxp=Z/p}
	  \ker(K_{x}^{*} \xrightarrow{\upsilon_{x}} \Z \xrightarrow{} \Zp)
	= \{ \lambda : \upsilon_{x}(\lambda) \equiv 0 \bmod p\} = K_{x}^{*p},
\end{equation}
induces an isomorphism $K_{x}^{*}/K_{x}^{*p} \simeq \Zp$. Now, we have the Kummer perfect pairing for $E_{x}$,
\begin{equation}
\label{E:Kummer pairing Kx}
	\pairing{g,\lambda}_{x} = \frac{g(\lambda^{1/p})}{\lambda^{1/p}} :
	\Gal(E_{x}/K_{x}) \times \faktor{K_{x}^{*}}{K_{x}^{*p}} \to \mu_{p},
\end{equation}
where $\lambda^{1/p}$ is any $p$th root of $\lambda$ in $E_{x}$.

Now, as we saw in \S\ref{subsec:covers as Galois ring extensions}, the isomorphism $\A_{Y} \simeq \A_{X} \otimes_{\Sigma} \Omega$ \eqref{E:cassels} is equivariant under the action of $\Gal(Y/X)$. Thus an automorphism $\tau$ of $Y/X$ can be restricted to each fiber of the cover $\pi : Y \to X$. In particular, for a ramification point $x \in \R$, $\tau$ induces an element $\tau_{x} \in \Gal(E_{x}/K_{x})$ in each local Galois group.

The explicit computations with local Kummer symbols carried out in~\cite[\S3.E]{adeles02} lead to the following algebraic definition of rotation data, valid over any algebraically closed base field $\k$ of characteristic prime to $p$.

\begin{definition}[Algebraic rotation data]
\label{D:algebraic rotation data}
The algebraic rotation data for an automorphism $\tau$ of a $p$-cyclic Galois cover $\pi : Y \to X$ with ramification locus $\R$ is defined by
\begin{equation}
\label{E:rotation data}
	\bigl(\log_{\zeta}\pairing{\tau_{x},z_{x}}_{x}\bigr)_{x \in \R} \in \prod_{x \in \R} \Zp,
\end{equation}
where $\log_{\zeta}$ is the discrete logarithm associated to a fixed primitive $p$th root of unity $\zeta$.
\end{definition}

\begin{proposition}
\label{P:rotation numbers}
The algebraic rotation data for a nontrivial automorphism $\tau$ of the superelliptic curve~\eqref{E:superelliptic} coincide with the classical analytic rotation data, under the assumption that $\chi(\tau) = \zeta$, where we identify $\cc_{p}$ with $\Gal(Y/X)$ as above.
\end{proposition}

\begin{proof}
This follows from~\cite[Proposition 3.79]{adeles02}, which together with our choices of $\chi$, $\tau$ and $\zeta$ show that the numbers $\log_{\zeta}\pairing{\tau_{x},z_{x}}_{x}$ are indeed the inverses modulo $p$ of the valuations $v_{x} = \upsilon_{x}(\f)$ and thus by \eqref{E:rotation is inverse of valuation} are equal to the rotation numbers.
\end{proof}

\begin{remark}
\label{R:rotation data}
Note that \eqref{E:rotation data} depends not only on the automorphism $\tau$ but also on the choice of primitive $p$th root of unity $\zeta$, and that the relation $\chi(\tau) = \zeta$ is needed for the equivalence in Proposition~\ref{P:rotation numbers}. In general it would hold only modulo a constant multiple in $\Zps$.

Thus, the equivalence class of $\R$-tuples as in \eqref{E:rotation data} modulo the action of $\Zps$ by multiplication (see for example~\cite[Theorem 3.71]{adeles02}) is independent of these choices.
\end{remark}



\end{document}